\documentclass[reqno,11pt]{amsart}

\usepackage{amsthm, amsmath, amssymb, stmaryrd}
\usepackage[utf8]{inputenc}
\usepackage[T1]{fontenc}

\usepackage[hyphenbreaks]{breakurl}
\usepackage[hyphens]{url}
\usepackage{systeme}
\usepackage[shortlabels]{enumitem}
\usepackage[hidelinks]{hyperref}
\usepackage{microtype}

\usepackage{bm}
\usepackage[margin=1in]{geometry}

\usepackage[textsize=scriptsize,backgroundcolor=orange!5]{todonotes}

\usepackage[noabbrev,capitalize,sort]{cleveref}
\crefname{equation}{}{}
\numberwithin{equation}{section}

\usepackage{mathtools}
\newtheorem{theorem}{Theorem}[section]
\newtheorem{proposition}[theorem]{Proposition}
\newtheorem{lemma}[theorem]{Lemma}

\theoremstyle{definition}
\newtheorem{definition}[theorem]{Definition}

\newtheorem{example}[theorem]{Example}

\theoremstyle{remark}
\newtheorem*{remark}{Remark}

\newcommand{\abs}[1]{\left\lvert#1\right\rvert}

\newcommand{\paren}[1]{\left( #1 \right)}

\DeclareMathOperator{\Span}{span}
\DeclareMathOperator{\mult}{mult}   
     
\DeclareMathOperator{\supp}{supp}

\newcommand*{\eqdef}{\stackrel{\mbox{\normalfont\tiny def}}{=}}

\newcommand{\EE}{\mathbb{E}}
\newcommand{\FF}{\mathbb{F}}

\newcommand{\RR}{\mathbb{R}}
\newcommand{\NN}{\mathbb{N}}
\newcommand{\ZZ}{\mathbb{Z}}
\newcommand*{\PP}{\mathbb{P}} 

\newcommand{\cL}{\mathcal L}
\newcommand{\cJ}{\mathcal J}
\newcommand{\cC}{\mathcal C}

\newcommand{\cF}{\mathcal F}

\newcommand{\cG}{\mathcal G}
\newcommand{\cM}{\mathcal M}

\newcommand{\cT}{\mathcal{T}}

\newcommand{\cH}{\mathcal{H}}

\newcommand{\Hasse}{\mathsf H} 

\newcommand{\totalB}{\mathbb B}
\newcommand{\prioB}{\mathcal B}

\newlength{\hght}

\makeatletter
\newcommand\thankssymb[1]{\textsuperscript{\@fnsymbol{#1}}}
\makeatother

\author[Ting-Wei Chao]{Ting-Wei Chao\thankssymb{1}}
\author[Hung-Hsun Hans Yu]{Hung-Hsun Hans Yu\thankssymb{2}}

\thanks{\thankssymb{1}Department of Mathematics, Massachusetts Institute of Technology, Cambridge, MA, USA. Email: {\tt twchao@mit.edu}}

\thanks{\thankssymb{2}Department of Mathematics, Princeton University, Princeton, NJ 08544\@.  Email: {\tt hansonyu@princeton.edu}}

\title{When Joints Meet Extremal Graph Theory: \linebreak Hypergraph Joints}

\begin{document}

\maketitle

\begin{abstract}
The Kruskal--Katona theorem determines the maximum number of $d$-cliques in an $n$-edge $(d-1)$-uniform hypergraph.
A generalization of the theorem was proposed by Bollob\'as and Eccles, called the partial shadow problem. 
The problem asks to determine the maximum number of $r$-sets of vertices that contain at least $d$ edges in an $n$-edge $(d-1)$-uniform hypergraph. 
In our previous work, we obtained an asymptotically tight upper bound via its connection to the joints problem, a problem in incidence geometry. 
In a different direction, Friedgut and Kahn generalized the Kruskal--Katona theorem by determining the maximum number of copies of any fixed hypergraph in an $n$-edge hypergraph, up to a multiplicative factor.

In this paper, using the connection to the joints problem again, we generalize our previous work to show an analogous partial shadow phenomenon for any hypergraph, generalizing Friedgut and Kahn's result.
The key idea is to encode the graph-theoretic problem with new kinds of joints that we call hypergraph joints.
Our main theorem is a generalization of the joints theorem that upper bounds the number of hypergraph joints, which the partial shadow phenomenon immediately follows from.
In addition, with an appropriate notion of multiplicities, our theorem also generalizes a generalization of H\"older's inequality considered by Finner.
\end{abstract}

\section{Introduction}
\subsection{The Kruskal--Katona theorem and its generalizations}
In extremal graph theory, we often try to understand quantitative relations between number of vertices, number of edges, and number of different subgraphs under various constraints.
One of the easiest examples asks the following: what is the maximum number of triangles in an $n$-edge graph?
In fact, the answer to this problem has long been known as a special case of the Kruskal--Katona theorem, which was independently proved by Kruskal \cite{Kru63} and Katona \cite{Kat68}.
Both of them actually solved the problem in more generality: in terms of hypergraphs, they determined the maximum number of $K_d^{(d-1)}$'s in a $(d-1)$-uniform hypergraph with $n$ edges for each $n$ and $d\geq 2$.
Despite the terminology we are using here, Kruskal and Katona were actually working in slightly different settings: Kruskal was motivated by counting faces of simplices, and Katona was studying set systems.
The Kruskal--Katona theorem has thus laid down important foundations for both algebraic combinatorics and the study of set systems.

Although the Kruskal--Katona theorem gives a very satisfying precise answer to the problem for all $n$, one slight drawback is the explicit expression of the bound is often difficult to work with.
As a remedy, Lov\'asz \cite{Lovasz93} observed that their bound can be weakened to get a bound that is still tight for infinitely many $n$'s and much easier to work with in general.
To state Lov\'asz's bound, for each positive integer $d$, we denote by $\binom{x}{d}$ the polynomial $\frac{x(x-1)\cdots (x-d+1)}{d!}$.

\begin{theorem}[Lov\'asz's version of the Kruskal--Katona theorem \cite{Lovasz93}]\label{thm:KK}
    For any real number $x\geq d$, any $(d-1)$-uniform hypergraph with $\binom{x}{d-1}$ edges has at most $\binom{x}{d}$ cliques on $d$ vertices.
\end{theorem}

Given the extremal-graph-theoretic motivation, it is natural to replace triangles or $K_d^{(d-1)}$ with general graphs or hypergraphs and ask the same question.
The $2$-uniform case was first considered by Alon \cite{Alon81}, and it was later extended to general hypergraphs (potentially with mixed uniformities) by Friedgut and Kahn \cite{FK98}.
To state their theorem, recall that the \emph{fractional edge-covering number} $\rho^*(\cH)$ of a hypergraph $\cH=(V,E)$ is the minimum possible value of $\sum_{e\in E}w(e)$ where $w$ is an assignment of nonnegative weights to the edges so that $\sum_{e\ni v}w(e)\geq 1$ for each $v\in V$.
Using the entropy method, Friedgut and Kahn showed the following.

\begin{theorem}[Friedgut--Kahn \cite{FK98}]\label{thm:FK}
    Let $\cH$ be a simple hypergraph with no isolated vertices.
    Then any hypergraph with $n$ edges contains at most $O_{\cH}(n^{\rho^*(\cH)})$ copies of $\cH$.
    Moreover, the bound is optimal up to a multiplicative factor depending only on $\cH$.
\end{theorem}

The technical assumption that $\cH$ has no isolated vertices is needed in order for $\rho^*(\cH)$ to be finite.

\subsection{The partial shadow phenomenon}
In the setting of the Kruskal--Katona theorem, a $d$-set is counted if all $d$ of its $(d-1)$-subsets are present in the given hypergraph.
Bollob\'as and Eccles \cite{BE15} first investigated what the maximum should be if we instead count the number of $r$-sets with at least $d$ of its $(r-1)$-subsets present in the hypergraph, where $r\geq d$ is some other parameter.
As the $(r-1)$-subsets of a $r$-set is usually called its \emph{shadow} in the field of set systems, they call the problem the \emph{partial shadow problem}.
Bollob\'as and Eccles made observations that the partial shadow problem is still related to the Kruskal--Katona theorem, and they conjectured that Lov\'asz's version of the Kruskal--Katona theorem continues to hold.
In particular, they conjectured that for any real number $x\geq d$, given any $\binom{x}{d-1}$ distinct $(r-1)$-sets, there are at most $\binom{x}{d}$ many $r$-sets that contain at least $d$ of the $(r-1)$-sets \cite[Conjecture 3]{BE15}.

Note that in Bollob\'as and Eccles' conjecture, the bound is the same as \cref{thm:KK}, which is independent of $r$.
To make the connection more rigorous, for any non-negative integer $t$, let $\cC_t(K_d^{(d-1)})$ be the $(d+t-1)$-uniform hypergraph on $[d+t]$ with edge set $\{[d+t]\backslash \{i\}\mid i\in [d]\}.$
Then in the partial shadow problem, we are actually counting the number of $r$-sets whose corresponding induced subgraphs contain $\cC_{r-d}(K_d^{(d-1)})$.
Motivated by this, we make the following definition.
\begin{definition}\label{def:M}
    For any hypergraph $\cH$ on $r$ vertices, we denote by $M(\cH,n)$ the maximum number of sets of $r$ vertices in a hypergraph with $n$ edges whose corresponding induced subgraphs contain $\cH$.
\end{definition}
With this definition, the partial shadow problem asks to determine $M(\cC_t(K_d^{(d-1)}),n)$, whereas the Kruskal--Katona theorem determines its precise value when $t=0$.

It is easy to see that $M(\cC_t(K_d^{(d-1)}),n)\geq M(K_d^{(d-1)},n)$, although the inequality is sometimes strict: Bollob\'as and Eccles showed that this is strict when $(n,d,t) = (12, 4, 1)$.
Nonetheless, their conjecture says that Lov\'asz's bound continues to hold.
This had been open until our previous work \cite{CY23}, which, surprisingly, proved the conjecture by considering a related problem in incidence geometry.
\begin{theorem}[{\cite[Theorem 1.7]{CY23}, conjectured in \cite[Conjecture 3]{BE15}}]\label{thm:partial-shadow-original}
    For any integers $d\geq 2$, $t\geq 0$ and any real number $x\geq d$,
    \[M\left(\cC_t(K_d^{(d-1)}), \binom{x}{d-1}\right)\leq \binom{x}{d}.\]
\end{theorem}

The only partial result before \cite{CY23} that we are aware of is Fitch's result \cite{Fitch18}, where Fitch proved the conjecture for $d=3$ and also showed that $M(\cC_t(K_3),n)= M(K_3,n)$ for all sufficiently large $n$ for any $t\in\ZZ_{\geq 0}$.

As mentioned before, the upper bound in \cref{thm:partial-shadow-original} does not depend on $t$, and setting $t=0$ recovers \cref{thm:KK}. 
One might ask if it is possible to define $\cC_t(\cH)$ appropriately for other hypergraphs $\cH$ so that $M(\cC_t(\cH),n)$ has an upper bound independent of $t$ which becomes the upper bound in \cref{thm:FK} when $t=0$.
Indeed, it is possible with the following definition.
\begin{definition}
Let $\cH=([d],E)$ be a hypergraph.
Its \emph{$t$-th cone} $\cC_t(\cH)$ is defined to be the hypergraph on $[d+t]$ with edge set $\{e\cup \{d+1,\ldots, d+t\}\mid e\in E\}$ for every non-negative integer $t$.
\end{definition}
It is clear that this definition is compatible with the earlier definition for $\cH = K_d^{(d-1)}$.
One of the main results of this paper is to show that \cref{thm:partial-shadow-original} indeed generalizes to other hypergraphs with this definition.
We call it \emph{the partial shadow phenomenon} for general hypergraphs.

\begin{theorem}[The partial shadow phenomenon]\label{thm:partial}
    Let $\cH=([d],E)$ be a simple hypergraph with no isolated vertices, and let $t\geq 0$ be an integer. Then
    \[M(\cC_t(\cH),n)=O_{\cH}\left(n^{\rho^*(\cH)}\right).\]
    In particular, the upper bound does not depend on $t$.
\end{theorem}

As in the previous work \cite{CY23}, the proof of \cref{thm:partial} also goes through some incidence geometry problem that we will introduce soon.
We remark that the proof of Friedgut and Kahn gives upper bounds with implicit constants growing roughly factorially in $t$.
As such, the only currently known proofs of Bollob\'as and Eccles' conjecture and the partial shadow phenomenon both go through incidence geometry. 
Since these problems are purely extremal-graph-theoretic, it is interesting to study whether there are proofs of these theorems without going through incidence geometry.

\subsection{The joints problem and its connection to the partial shadow problem}
As mentioned above, the proofs of both Bollob\'as and Eccles' conjecture and \cref{thm:partial} go through a problem in incidence geometry, which is in fact called the \emph{joints problem}.
Here, we give a brief introduction of the joints problem, and explain how it is related to the partial shadow problem.

The joints problem asks for the maximum number of joints formed by $n$ lines in $\FF^d$, where a joint is a point contained in $d$ lines with linearly independent directions. 
The problem first appeared in \cite{CEGPSSS92}, and attention was later brought to this problem by Wolff \cite{Wol99} due to its connection to the Kakeya problem in harmonic analysis.
Despite its simple formulation, it was already difficult for the special case $d=3$ and $\FF=\RR$, which was only resolved a decade later by Guth and Katz \cite{GK10} up to a constant multiplicative factor using the polynomial method introduced by Dvir \cite{Dvir09}. 
Soon after their work, it was quickly generalized by others \cite{ CI14, Dvir10, KSS10, Qui09, Tao14}, showing that $n$ lines in $\FF^d$ can determine $O_d(n^{\frac{d}{d-1}})$ joints. 
On the other hand, for any positive integer $m\geq d$, there is a construction with $\binom{m}{d}$ joints using $\binom{m}{d-1}$ lines as follows.
Consider $m$ hyperplanes in general position. Take our set of lines to be the collection of all $(d-1)$-wise intersections of these hyperplanes. 
It follows that any $d$-wise intersection of these hyperplanes is a joint. 
Therefore, the construction matches the upper bound obtained from the polynomial method, up to a multiplicative factor. 
There have been recent works that aimed to close up the gap between the upper bound and the lower bound.
Yu and Zhao \cite{YZ23} first proved an upper bound that is tight up to a $(1+o(1))$-factor, and later, Chao and Yu \cite{CY23} tightened up the upper bound even more, showing that $\binom{x}{d-1}$ lines in $\FF^d$ form at most $\binom{x}{d}$ joints for any real number $x\geq d$.

The connection of the joints problem to the Kruskal--Katona theorem was first made explicit by Yu and Zhao \cite{YZ23}, where they introduced \emph{generically induced configurations}.
A generically induced configuration is a special configuration that can be formed by taking hyperplanes $H_1,\dots,H_m$ in general position, and using only points and lines of the form $\cap_{i\in I}H_i$ for some $I\subseteq [m]$. 
If we restrict ourselves to only generically induced configurations, the original joints problem in a $3$-dimensional space becomes determining $M(K_3,n)$.
To see this, note that each edge $\{i,j\}\in \binom{[m]}{2}$ corresponds to two hyperplanes $H_i,H_j$, which intersect in a line $\ell_{ij}$.
Any triangle $\{i,j,k\}$ corresponds to a point $H_i\cap H_j\cap H_k$ in the same way, and it is clear that it is a joint formed by $\ell_{ij},\ell_{jk}$ and $\ell_{ik}$.
Therefore, bounds for the joints problem translate to bounds for the Kruskal--Katona theorem.
In particular, the bound in our previous work \cite{CY23} translates to Lov\'asz's version of the Kruskal--Katona theorem.

In fact, we may also consider configurations of a more general form called \emph{projected generically induced configurations}, which was introduced in our previous work \cite{CY23}. A projected generically induced configuration in $\FF^d$ is a configuration obtained by taking a generically induced configuration in some higher dimensional space $\FF^{d+t}$, and projected back to $\FF^d$ via a generic projection.
It turns out that through a similar correspondence, when applied to projected generically induced configurations, our previous result \cite{CY23} proves Bollob\'as and Eccles' conjecture \cite[Conjecture 3]{BE15}.

\subsection{Other variants of the joints problem}
Given the correspondence above via generically induced configurations, it is natural to ask whether there is a corresponding incidence geometry problem to Friedgut and Kahn's result through the same correspondence.
The corresponding geometric configuration turns out to have flats, i.e. affine subspaces, of higher dimension, as the hypergrpah $\cH=([d],E)$ might have edges with uniformity $d-2$ or lower.
Luckily, joints involving objects of higher dimensions have been considered, which we briefly introduce here.

To generalize, a point $p$ is a joint of flats $F_1,\dots,F_r$ in $\FF^d$ if $p$ is their common intersection, $d=\dim F_1+\cdots+\dim F_r$ and they do not lie in any hyperplane simultaneously. 
Equivalently, there exists a linear transformation that sends the origin to $p$ and $\Span\{\vec{e_{j}}\mid j\in I_i\}$ to $F_i$ for some $I_1,\dots,I_r$ that partition $[d]$.
The first major progress on this problem was made by Yang \cite{Yangthesis} in his thesis, who determined the maximum number of joints of flats over $\RR$ up to $o(1)$-losses in the exponents.
The problem was later solved by Tidor, Yu, and Zhao \cite{TYZ22}, who determined the correct order of magnitude for these generalizations and, more generally, for joints of varieties.

There are other variants of the joints problem that have also been considered, including a colorful version that is more often referred to as the \emph{multijoints problem}. 
Suppose that we are given a set of lines $\cL_i$ with color $i$ for each $i\in[d]$, the multijoints problem asks: what is the maximum number of \emph{rainbow joints}, i.e. joints formed by $d$ lines with distinct colors, in terms of $\abs{\cL_1},\dots,\abs{\cL_d}$? 
Zhang \cite{Zha20} showed that the number is at most $O\left((\abs{\cL_1}\dots\abs{\cL_d})^{\frac{1}{d-1}}\right)$, and this bound is tight up to a multiplicative factor.
The upper bound was also improved by Yu and Zhao \cite{YZ23}, and later Chao and Yu \cite{CY23} showed that the upper bound is tight up to a multiplicative $(1+o(1))$-factor as long as each joint is counted with an appropriately defined multiplicity, which we will discuss later as well.
The corresponding multijoints problem for joints of flats (or varieties in general) was also resolved by Tidor, Yu and Zhao \cite{TYZ22}, whose bound is also tight up to a multiplicative $(1+o(1))$-factor if the joints are counted with appropriate multiplicities.

\subsection{Encoding Friedgut--Kahn using joints}
In this subsection, we will state our main theorem that encodes the Friedgut--Kahn theorem using joints, which is the key to prove the partial shadow phenomenon (\cref{thm:partial}).
Before stating the theorem, we first need to introduce the setting of the hypergraph $\cH$ we will work on and some notions on it.

\begin{definition}
    Let $\cH=([d],E,c)$ be an edge-colored (multi-)hypergraph with vertex set $[d]$, edge set $E$ (which may be a multi-set), and an edge coloring $c:E\rightarrow [r]$. For integers $1\leq k_1,\dots,k_r\leq d-1$, we say that the coloring $c$ is \emph{$(k_1,\dots,k_r)$-uniform} if the preimage $c^{-1}(i)$ is a non-empty $k_i$-uniform simple hypergraph for each color $i\in [r]$. In this case, we say $\cH$ is \emph{uniformly edge-colored}.

    For a weight function $w:E\rightarrow \RR_{\geq 0}$, we say that $w$ \emph{covers} $\cH$ if 
    \[\sum_{e\ni j}w(e)\geq 1\quad\forall j\in [d].\]
    Denote by $\abs{w}$ the sum $\sum_{e\in E}w(e)$. Also, for each color $i\in [r]$, set
    \[\bar{w}_i=\sum_{e\in c^{-1}(i)}w(e),\]
    and we say that the sequence $(\bar{w}_1,\dots,\bar{w}_r)$ is the \emph{subtotal sequence} of $w$.

\end{definition}

Given the information of $\cH$, we may define an $\cH$-joint. We shall keep in mind that each vertex of the graph corresponds to a hyperplane, and hence an edge of uniformity $d-k$ corresponds to a $k$-flat.

\begin{definition}
    Let $\cH=([d],E,c)$ be an edge-colored (multi-)hypergraph where the coloring $c$ is $(d-k_1,\dots,d-k_r)$-uniform and $1\leq  k_1,\dots,k_r\leq d-1$. For each $i$, let $\cF_i$ be a multi-set of $k_i$-flats (i.e. $k_i$-dimensional affine subspaces) in $\FF^d$. We say that a point $p\in \FF^d$ is an \emph{$\cH$-joint} if there exist flats $(F_{p,e})_{e\in E}$ with $F_{p,e}\in \cF_{c(e)}$ and an invertible affine transformation $A:\FF^d\to \FF^d$ such that $A(0)=p$ and $F_{p,e}=A(\Span\{\vec{e_j}\mid j\notin e\})$, where $\vec{e_j}$ is the $j$-th unit vector. In this case, we say that $A$ \emph{witnesses the $\cH$-joint $p$} and the flats $\left(F_{p,e}\right)_{e\in E}$ \emph{form the $\cH$-joint $p$}.

    Let $\cJ$ be the set (without multiplicities) of $\cH$-joints. We say that $(\cJ,\cF_1,\dots,\cF_r)$ is an \emph{$\cH$-joints configuration of dimension $(k_1,\dots,k_r)$}.
\end{definition}

We provide several examples here to make it easier to digest the definition.

\begin{example}
    Let $d\geq 2$ be a positive integer, and let $\cH = K^{(d-1)}_d$.
    If all edges in $\cH$ are of the same color, then $\cH$ is $(d-1)$-uniform, showing that we should take $k_1=1$.
    Therefore, for a multi-set $\cL$ of $1$-flats (i.e. lines) in $\FF^d$, a point $p$ is an $\cH$-joint formed by the lines $(\ell_{[d]\backslash i})_{i\in [d]}$ if there is an invertible affine transformation $A:\FF^d\to \FF^d$ sending $0$ to $p$ and $\textup{span}(\vec{e_i})$ to $\ell_{[d]\backslash i}$.
    It is then clear that $(\ell_{[d]\backslash i})_{i\in [d]}$ form an $\cH$-joint at $p$ if and only if they all pass through $p$ and they have linearly independent directions.
    This matches the usual definition of joints.

    If $\cH$ is colored so that $c([d]\backslash i)=i$ for each $i\in[d]$ instead, then $\cH$ is $(d-1,\ldots, d-1)$-uniform, showing that we should take $k_1=\cdots=k_d=1$.
    In this case, for multi-sets $\cL_1,\ldots, \cL_d$ of lines, a point $p$ is an $\cH$-joint formed by the lines $(\ell_{[d]\backslash i})_{i\in [d]}$ if they form a joint there in the usual sense and also  $\ell_{[d]\backslash i}\in \cL_i$ for each $i\in[d]$.
    Therefore an $\cH$-joint in this case is the same as a multijoint of lines.
\end{example}

\begin{example}\label{ex:joints-of-flats}
    Let $\cH$ be the $4$-uniform hypergraph on $[6]$ with edges $\{1,2,3,4\}, \{1,2,5,6\}$ and $\{3,4,5,6\}$ with the same color.
    Then $k_1=2$, and for a multi-set $\cF$ of $2$-flats in $\FF^6$, a point $p$ is an $\cH$-joint formed by $2$-flats $F_1$, $F_2$, $F_3$ in $\cF$ if the following holds: there exists an invertible affine transformation $A:\FF^6\to\FF^6$ sending $0$ to $p$, $\textup{span}\{\vec{e_1}, \vec{e_2}\}$ to $F_1$, $\textup{span}\{\vec{e_3},\vec{e_4}\}$ to $F_2$ and $\textup{span}\{\vec{e_5},\vec{e_6}\}$ to $F_3$.
    Note that this is equivalent to that $F_1,F_2,F_3$ do not lie in the same hyperplane, so an $\cH$-joint is the same as a joint of $2$-flats in $\FF^6$.

    In general, an $\cH$-joint is a joint of flats if $\cH=([d],E)$ where $\bigcup_{e\in E}([d]\backslash e)=[d]$ is a disjoint union. 
\end{example}

\begin{example}
    Let $\cH=([5],E)$ be the $5$-cycle with $E=\{\{1,2\},\{2,3\},\{3,4\},\{4,5\},\{1,5\}\}$, and all edges of $\cH$ have the same color.
    Then $k_1=5-2=3$ as $\cH$ is $2$-uniform.
    For a multi-set $\cF$ of $3$-flats in $\FF^5$, a point $p$ is an $\cH$-joint formed by $3$-flats $F_1,F_2,F_3,F_4,F_5$ in $\cF$ if there is an invertible transformation $A:\FF^5\to \FF^5$ satisfying the following: $A$ sends $0$ to $p$ and $\textup{span}(\vec{e_{i+2}},\vec{e_{i+3}}, \vec{e_{i+4}})$ to $F_i$ for each $i\in[5]$, where the indices are interpreted modulo $5$.
    Note that in this case, we know that for each $i\in[5]$, $F_i\cap F_{i+2}$ is exactly $A(\textup{span}(\vec{e_{i+4}}))$, which also lies in $F_{i+1}$.
    Therefore an equivalent formulation is that $F_1,\ldots, F_5$ are $3$-flats containing $p$ such that $(F_i\cap F_{i+2})_{i\in [5]}$ are linearly independent lines, and $F_i\cap F_{i+2}\subseteq F_{i+1}$ for each $i\in[5]$.

    We remark that the equivalent description is specific to the $5$-cycles and is worked out only to help the readers understand the definition in this particular case.
\end{example}

Now, we may state our $\cH$-joints theorems that encode the Friedgut--Kahn theorem. We start with simple joints, i.e. counting each joint with multiplicity $1$.

\begin{theorem}[Simple $\cH$-joints]\label{thm:HSimpleJoints}
    Let $(\cJ,\cF_1,\dots,\cF_r)$ be an $\cH$-joints configuration of dimension $(k_1,\dots,k_r)$, where $\cH=([d],E,c)$ is a uniformly edge-colored (multi-)hypergraph. Fix a weight function $w:E\rightarrow\RR_{\geq 0}$ that covers $\cH$ and let $(\bar{w}_1,\dots,\bar{w}_r)$ be its subtotal sequence. Then
    \[\abs{\cJ}\leq C_{\cH,w}\abs{\cF_1}^{\bar{w}_1}\dots\abs{\cF_r}^{\bar{w}_r},\]
    where
    \[C_{\cH,w}=d!^{\abs{w}-1}\prod_{i=1}^r \left(\frac{1}{k_i!}\right)^{\bar{w}_i}\prod_{e\in c^{-1}(i)}\left(\frac{w(e)}{\bar{w}_i}\right)^{w(e)}.\]
    Here, we are using the convention that $\left(\frac{w(e)}{\bar{w}_i}\right)^{w(e)}=1$ when $w(e)=0$, even if $\bar{w}_i=0$.
\end{theorem}

To immediately see how this is related to the bound in Friedgut and Kahn's result (\cref{thm:FK}), note that the total sum of the exponents is $\abs{w}$, whose minimum is exactly the fractional covering number $\rho^*(\cH)$.
We will utilize this to prove \cref{thm:partial} as well.

The constant in \cref{thm:HSimpleJoints} might not be tight, but similar to \cite{CY23}, the constant would be tight in many cases once we count the joints with appropriate multiplicities.
The definition of multiplicity might seem technical and artificial at first glance as it is an artifact of our proof technique, but we will show connections to other quantities later.

\begin{definition}\label{def:eta}
    In the same setup as \cref{thm:HSimpleJoints}, for any joint $p$, let $\cT_p$ be the set of tuples $(F_{p,e})_{e\in E}$ with $F_{p,e}\in \cF_{c(e)}$ that form an $\cH$-joint at $p$. We define the multiplicity $\eta(p)$ as follows.  
    For each $i\in [r]$, let $e_i$ be an edge randomly and independently chosen from $c^{-1}(i)$ so that the probability $\PP(e_i=e) = w(e)/\bar{w}_i$ for each $e\in c^{-1}(i)$. We set $\eta(p)$ to be the non-negative real number such that
    \[\log_2\eta(p)=\max \sum_{i=1}^r\bar{w_i}\left(H(F_{p,e_i})-H(e_i)\right),\]
    where $H(F_{p,e_i})$ and $H(e_i)$ denote the Shannon entropy (see \cref{section:tool}) of the random variables $F_{p,e_i}$ and $e_i$, respectively, and the maximum is taken over all possible distributions of random tuple of flats $(F_{p,e})_{e\in E}\in\cT_p$ independent from $(e_i)_{i\in[r]}$. 
    Note that the maximum can be achieved as the space of distributions is compact.
\end{definition}

\begin{remark}
    This definition of multiplicity is actually closely related to the one introduced in \cite{CY23}. 
    To see this, recall from \cref{ex:joints-of-flats} that when $\cH=([d],E,c)$ where the complements of edges are disjoint and the coloring $c$ gives each uniformity a color, an $\cH$-joint is just a (multi-)joint of flats. 
    Under the notation used in \cite[Definition 6.4]{CY23}, if we set the weight on each edge $w(e)=\frac{1}{m_1+\dots+m_r-1}$, the multiplicity $\eta(p)$ we defined equals $\nu^*(p)^{(m_1+\dots+m_r)/(m_1+\dots+m_r-1)}$.

    The multiplicity is also chosen so that if $p$ is a simple joint, i.e. if $\cT_p$ has size $1$, then $\eta(p) = 1$.
\end{remark}

By counting each joint $p$ with multiplicity $\eta(p)$, we can show that not only the same upper bound holds, but also the constant $C_{\cH,w}$ cannot be improved for many cases.

\begin{theorem}[$\cH$-joints with multiplicities]\label{thm:HMultiplicityJoints}
    Let $(\cJ,\cF_1,\dots,\cF_r)$ be an $\cH$-joints configuration of dimension $(k_1,\dots,k_r)$, where $\cH=([d],E,c)$ is a uniformly edge-colored (multi-)hypergraph. Fix a weight function $w:E\rightarrow\RR_{\geq 0}$ that covers $\cH$ and let $(\bar{w}_1,\dots,\bar{w}_r)$ be its subtotal sequence. Then
    \[\sum_{p\in\cJ}\eta(p)\leq C_{\cH,w}\abs{\cF_1}^{\bar{w}_1}\dots\abs{\cF_r}^{\bar{w}_r}\]
    holds for the same constant $C_{\cH,w}$ as in \cref{thm:HSimpleJoints}.

    Moreover, the constant $C_{\cH,w}$ is the best possible constant if the weight $w$ satisfies $\sum_{e\ni j}w(e)= 1$ for all $j\in [d]$. 
\end{theorem}

\subsection{Proof strategy}
Our proof builds on Tidor, Yu and Zhao's proof of the joints of varieties theorem \cite{TYZ22}.
Similar to the proof there, we will assign appropriate vanishing conditions for joints on each flat according to a ``handicap'' that will be assigned later.
As we will apply properties from \cite{TYZ22} as black boxes for this part, we refer interested readers to \cite[Section 2]{TYZ22} for a more detailed discussion of the ideas there.
The main new idea in our paper is a more general way of combining the vanishing conditions from different flats, and a new vanishing lemma that is tailored to the hypergraph joints.
This is done in \cref{subsec:LW-step}.

An additional idea is a different way to end the proof.
In the proof in \cite{TYZ22}, after the polynomial method was done to obtain a key inequality, the proof was concluded with a series of applications of the weighted AM-GM inequality.
Here, we replace the technical deduction with a more conceptual one, showing that the key inequality actually implies the following entropic inequality.

\begin{theorem}\label{thm:entropy-key-ineq}
    Let $\cH$, $(\cJ, \cF_1,\ldots, \cF_r)$ and $(\cT_p)_{p\in\cJ}$ be the same as in \cref{thm:HSimpleJoints} and \cref{def:eta}.    
    Let $p\in \cJ$ be a randomly chosen $\cH$-joint and let $(F_{p,e})_{e\in E}\in\cT_p$ be randomly chosen. 
    For each $i\in [r]$, let $e_i$ be an edge randomly and independently chosen from $c^{-1}(i)$ so that $\PP(e_i=e) = w(e)/\bar{w}_i$ for each $e\in c^{-1}(i)$.
    Then
    \[H(p)+\sum_{i=1}^r\bar{w}_i\left(H(F_{p,e_i}\mid p)-H(e_i)\right)\leq \sum_{i=1}^r\bar{w}_iH(F_{p,e_i})+\log_2 C_{\cH,w}\]
    where $C_{\cH,w}$ is the same as in \cref{thm:HSimpleJoints}.
\end{theorem}

It turns out that when restricted to certain configurations, this theorem would almost become Shearer's inequality \cite{CGFS86} (with the only difference being an additive constant).
Therefore we nickname this theorem the \emph{geometric Shearer's inequality}, and we will make this connection more precise in \cref{subsec:shearer}.
Once this is proven, \cref{thm:HSimpleJoints} and \cref{thm:HMultiplicityJoints} can be easily seen as immediate corollaries of the geometric Shearer's inequality.
\cref{thm:partial} will then follow from \cref{thm:HSimpleJoints}.

It is surprising that the key inequality obtained via the polynomial method would naturally imply some entropic inequality, and this might not be a coincidence as Friedgut and Kahn also used the entropy method to prove their result.
This is hinting that there might be a deeper connection between the polynomial method and the entropy method, and because of this, we believe that the geometric Shearer's inequality may be of independent interest.

\subsection{Axis-parallel configurations}
Previously, we drew connections from the joints problem to extremal graph theory using generically induced configurations and projected generically induced configurations.
In this part, we will draw more connections using more specific configurations, which we call the \emph{axis-parallel configurations}.
For any hypergraph $\cH$ equipped with an edge-coloring $c$ such that each edge receives a unique color, axis-parallel configurations are configurations of the form $(\cJ,\cF_1,\dots,\cF_r)$ where $\cF_{i}$ (which may be a multi-set) consists of flats parallel to $\Span (\vec{e_j}\mid j\notin e)$ if $e$ is the edge with color $i$. When restricted to axis-parallel configurations, \cref{thm:HMultiplicityJoints} becomes the discrete version of a generalization of H\"older's inequality proved by Finner \cite{Fin92} up to a multiplicative factor.

\begin{theorem}[Generalized H\"older's inequality \cite{Fin92}]\label{thm:GenHolder}
    Let $I_1,\dots,I_m\subsetneq [d]$ be non-trivial subsets. Suppose $w_i> 0$ are weights such that $\sum_{I_i\ni j}w_i\geq 1$ for all $j$ in $[d]$.
    Let $S$ be a finite set.
    For every $i\in[m]$, let $\pi_i:S^d\rightarrow S^{I_i}$ be the projection onto the entries indexed by $I_i$. 
    Also, let $f_i:S^{I_i}\rightarrow\RR_{\geq 0}$ be some functions.
    Then
    \[\sum_{p\in S^d}\left(\prod_{i=1}^m f_i(\pi_i(p))^{w_i}\right)\leq \prod_{i=1}^m\left(\sum_{p_i\in S^{I_i}}f_i(p_i)\right)^{w_i}.\]
\end{theorem}

In \cref{section:GenHolder}, we will show this inequality using \cref{thm:HMultiplicityJoints} and a tensor power trick.

In addition, when restricted to axis-parallel configurations, \cref{thm:entropy-key-ineq} becomes Shearer's inequality up to an additive constant.
This justifies the nickname ``geometric Shearer's inequality''.

\begin{theorem}[Shearer's inequality \cite{CGFS86}]\label{thm:Shearer}
Let $I_1,\dots,I_m\subsetneq [d]$ be non-trivial subsets. Suppose $w_i> 0$ are weights such that $\sum_{I_i\ni j}w_i\geq 1$ for all $j$ in $[d]$.
Suppose $X_1,\dots,X_d$ are random variables with finite support. Denote by $X_I$ the tuple $(X_i)_{i\in I}$. Then we have
\[H(X_1,\dots,X_d)\leq \sum_{i=1}^m w_iH(X_{I_i}).\]
\end{theorem}

Again, we will show Shearer's inequality using the geometric Shearer's inequality and another tensor power trick in \cref{section:GenHolder}.

\subsection*{Paper organization}
In \cref{section:tool}, we give a brief review of the Shannon entropy as well as some useful inequalities. In \cref{section:keyineq}, we state a key inequality, \cref{thm:key-ineq}, and prove \cref{thm:entropy-key-ineq,thm:HSimpleJoints,thm:HMultiplicityJoints} assuming the key inequality. In \cref{section:poly}, we use the polynomial method to prove the key inequality. In \cref{section:graphhom,section:GenHolder}, we focus on the applications of our $\cH$-joints theorems on the number of subgraphs and inequalities in analysis, respectively.

\section{Toolbox}\label{section:tool}
In this section, we introduce the tools we need throughout the paper. We will cover some basics about Shannon entropy, and also some useful inequalities as well. For a more detailed introduction to Shannon entropy, see \cite[Section 14.6]{AS00}.

\begin{definition}[Shannon entropy]
    Let $X$ be a random variable such that the support $\supp(X)$ is finite. The \emph{entropy of $X$} is given by
    \[H(X)\eqdef \quad\smashoperator{\sum_{x\in \supp(X)}}\ -\PP(X=x)\log_2 \PP(X=x).\]
    For random variables $X_1,\dots,X_d$ with finite support, the entropy $H(X_1,\dots,X_d)$ is defined to be the entropy of the random tuple $(X_1,\dots,X_d)$.
\end{definition}

We will also need the definition of conditional entropy. There are several equivalent definitions.

\begin{definition}[conditional entropy]
    Let $X,Y$ be random variables with finite support. The \emph{conditional entropy of $X$ given $Y$} is defined by
    \begin{align*}
        H(X\mid Y)\eqdef &H(X,Y)-H(Y)\\
        =&\quad\smashoperator{\sum_{y\in\supp(Y)}}\ \PP(Y=y)H(X\mid Y=y)\\
        =&\quad\smashoperator{\sum_{(x,y)\in \supp(X,Y)}}\ -\PP(X=x,Y=y)\log_2 \frac{\PP(X=x,Y=y)}{\PP(Y=y)},
    \end{align*}
    where $H(X\mid Y=y)$ is the entropy of the random variable $X\mid Y=y$.
\end{definition}

The following uniform bound will be useful when we upper bound some entropy terms.
\begin{proposition}[uniform bound]
    Let $X$ be a random variable with finite support. We have
    \[H(X)\leq \log_2\abs{\supp(X)}.\]
    The equality holds if and only if $X$ is uniform.
\end{proposition}

The following lemma gives an alternative way to upper bound entropies. 
This is how we are going to connect a statement proven by the polynomial method to the geometric Shearer's inequality.
We provide a proof as this is not commonly used in the literature.

\begin{lemma}\label{lem:jensen}
    Let $S,T$ be finite sets, and let $(x_{s,t})_{s\in S,t\in T}$ be non-negative reals.
    Suppose that $A$ is a positive real such that $\sum_{t\in T} x_{s,t}\leq A$ for all $s\in S$.
    Then for any joint random variables $(\Tilde{s}$, $\Tilde{t})$ taking values in $S\times T$,
    \[\EE[-\log_2 x_{\Tilde{s},\Tilde{t}}]\geq H(\Tilde{t}\mid \Tilde{s})-\log_2A.\]
\end{lemma}
\begin{proof}
    We first fix an $s\in S$, and let $p_t$ be the probability that $\Tilde{t}=t$ conditioning on $\Tilde{s}=s$.
    Then by Jensen's inequality and the convexity of the function $-\log x$,
    \begin{align*}
        \EE[-\log_2x_{\Tilde{s},\Tilde{t}}\mid \Tilde{s}=s] =& -\sum_{t\in T,\,p_t>0}p_t\log_2 x_{s,t} = -\sum_{t\in T,\,p_t>0}p_t\log_2 \left(\frac{x_{s,t}}{p_t}\right)-\sum_{t\in T,\,p_t>0}p_t\log_2 p_t\\
        \geq& -\log_2\left(\sum_{t\in T,\,p_t>0}p_t\cdot \frac{x_{s,t}}{p_t}\right)+H(\Tilde{t}\mid \Tilde{s}=s)\\
        \geq& -\log_2 A+H(\Tilde{t}\mid \Tilde{s}=s).
    \end{align*}
    Therefore
    \[\EE[-\log_2x_{\Tilde{s},\Tilde{t}}] \geq \EE_{s\sim \Tilde{s}}H(\Tilde{t}\mid \Tilde{s}=s)-\log_2 A = H(\Tilde{t}\mid \Tilde{s})-\log_2 A.\qedhere\]
\end{proof}

We will also need the unweighted version of \cref{thm:GenHolder}, i.e. the discrete Loomis--Whitney inequality \cite{LW49}.
\begin{theorem}[Loomis--Whitney \cite{LW49}]
    Let $I_1,\dots,I_m\subsetneq [d]$ be non-trivial subsets. Suppose $w_i> 0$ are weights such that $\sum_{I_i\ni j}w_i\geq 1$ for every $j\in [d].$
    Let $S$ be a finite set.
    Let $\pi_i:S^d\rightarrow S^{I_i}$ be the projection onto the entries indexed by $I_i$. Let $T$ be a subset of $S^d$.
    Then
    \[\abs{T}\leq \prod_{i=1}^m\abs{\pi_i(T)}^{w_i}.\]
\end{theorem}

\section{Key inequality}\label{section:keyineq}
In this section, we state the key inequality and show how it implies the main theorems.
We defer the proof of the key inequality until the next section as we would like to present the part of the proof without polynomials first.
\begin{theorem}\label{thm:key-ineq}
    Let $(\cJ, \cF_1,\ldots, \cF_r)$ be any $\cH$-joints configuration.
    Let $w$ be a weight covering $\cH$.
    Then for any $W:\cJ\to \RR_{\geq 0}$ with $\sum_{p\in \cJ}W(p)=\frac{1}{d!}$, there exists $(b_{p,F})_{p\in  F}$ where $p\in \cJ$, $F\in \cF_1\cup\cdots\cup \cF_r$ and $b_{p,F}\geq 0$ such that the following holds.
    \begin{enumerate}
        \item  For any $p\in \cJ$ and $(F_{p,e})_{e\in E}\in \cT_p$,
        \[\left(\prod_{e\in E}b_{p,F_{p,e}}^{w(e)}\right)^{\frac{1}{\abs{w}-1}}\geq W(p).\]
        
        \item For any $F\in \cF_1\cup\cdots\cup \cF_r$,
        \[\sum_{p\in F\cap \cJ}b_{p,F}\leq \frac{1}{(\dim F)!}.\]
    \end{enumerate}
\end{theorem}

Assuming \cref{thm:key-ineq}, we give a proof of the geometric Shearer's inequality, \cref{thm:entropy-key-ineq}.
\begin{proof}[Proof of \cref{thm:entropy-key-ineq} assuming \cref{thm:key-ineq}]
    In this proof, we rewrite $F_{p,e_i}$ as $F_i$ for simplicity.
    By taking $W(p^*)=\frac{\PP(p=p^*)}{d!}$ and applying \cref{thm:key-ineq} to get the $b_{p,F}$'s, we have
     \begin{align*}
         H(p)+\log_2 d!=&\EE[-\log_2 W(p)]\\
         \geq&\EE\left[-\sum_{e\in E}\frac{w(e)}{\abs{w}-1}\log_2b_{p,F_{p,e}}\right]\\
         =&\sum_{i=1}^r\frac{\bar{w}_i}{\abs{w}-1}\EE\left[-\log_2b_{p,F_i}\right]\\
         \geq& \sum_{i=1}^r\frac{\bar{w}_i}{\abs{w}-1}\left(H(p\mid F_i)+\log_2k_i!\right)
     \end{align*}
    where the last inequality follows from \cref{lem:jensen}, where we set $b_{p^*,F}=0$ if $p^*\notin F$, and the fact that $\sum_{p^*\in \cJ\cap F_i}b_{p^*,F_i}\leq 1/k_i!$.

    Note that $H(p\mid F_i)=H(p,F_i)-H(F_i)=H(p)+H(F_i\mid p)-H(F_i)$. Thus, we get
    \[H(p)+\log_2 d!\geq \sum_{i=1}^r\frac{\bar{w}_i}{\abs{w}-1}\left(H(p)+H(F_i\mid p)-H(F_i)+\log_2k_i!\right).\]
    By multiplying $\abs{w}-1$ on both sides and rearranging, we get
    \begin{align*}
        &H(p)+\sum_{i=1}^r\bar{w}_i\left(H(F_i\mid p)-H(e_i)\right)\\
        \leq &\sum_{i=1}^r\bar{w}_iH(F_i)+(\abs{w}-1)\log_2 d!-\sum_{i=1}^r\bar{w}_i\log_2 k_i!-\sum_{i=1}^r \bar{w}_iH(e_i)\\
        =&\sum_{i=1}^r\bar{w}_iH(F_i)+\log_2 C_{\cH,w}.\qedhere
    \end{align*}
\end{proof}

As hinted in the introduction, the geometric Shearer's inequality immediately implies the two $\cH$-joints theorems.

\begin{proof}[Proof of \cref{thm:HSimpleJoints} assuming \cref{thm:entropy-key-ineq}]
    For each joint $p^*$, fix a deterministic choice of $(F_{p^*,e})_{e\in E}$ that forms an $\cH$-joint at $p^*$ beforehand. 
    Now we let $p$ be a uniformly chosen joint.
    Then we have $H(F_{p,e_i}\mid p)=H(e_i)$, and hence \cref{thm:entropy-key-ineq} gives
    \begin{align*}
        \log_2\abs{\cJ}=H(p)\leq\sum_{i=1}^r\bar{w}_iH(F_{p,e_i})+\log_2 C_{\cH,w}        \leq \log_2 \left(C_{\cH,w}\abs{\cF_1}^{\bar{w_1}}\dots\abs{\cF_r}^{\bar{w_r}}\right),
    \end{align*}
    where the second inequality follows from the uniform bound.
\end{proof}

\begin{proof}[Proof of \cref{thm:HMultiplicityJoints} assuming \cref{thm:entropy-key-ineq}]
    For each $p^*\in \cJ$, we let $(F_{p^*,e})_{e\in E}\in\cT_{p^*}$ be a random tuple such that $\log_2 \eta(p^*)=\sum_{i=1}^r\bar{w_i}\left(H(F_{p^*,e_i})-H(e_i)\right)$. We pick a random $p\in\cJ$ independently such that
    \[\PP(p=p^*)=\frac{\eta(p^*)}{\sum_{p'\in \cJ}\eta(p')}.\]
    By the definition of entropy and conditional entropy, we have
    \begin{align*}
        H(p)+\sum_{i=1}^r\bar{w_i}\left(H(F_{p,e_i}\mid p)-H(e_i)\right)=&\sum_{p^*\in\cJ}\left(-\PP(p=p^*)\log_2 \PP(p=p^*)+\PP(p=p^*)\log_2\eta(p^*)\right)\\
        =&\sum_{p^*\in\cJ}\left(\PP(p=p^*)\log_2\left(\sum_{p'\in \cJ}\eta(p')\right)\right)\\
        =&\log_2\left(\sum_{p'\in \cJ}\eta(p')\right).
    \end{align*}
    Thus, \cref{thm:entropy-key-ineq} gives
    \[\log_2\left(\sum_{p'\in \cJ}\eta(p')\right)\leq \sum_{e\in E}w(e)H(F_{p,e})+\log_2 C_{\cH,w}\leq \log_2 \left(C_{\cH,w}\abs{\cF_1}^{\bar{w_1}}\dots\abs{\cF_r}^{\bar{w_r}}\right),\]
    where the second inequality follows from the uniform bound.
    
    Next, we show that the constant is the best possible when $\sum_{e\ni j}w(e)= 1$ for all $j\in [d]$.
    This can be seen by considering the configuration as follows. Let $m$ be a positive integer that will be taken to tend to infinity. Let $H_1,\dots,H_m$ be hyperplanes in $\FF^d$ in general position. Take $\cJ$ to be the set of all the $d$-wise intersection of the hyperplanes, and take $\cF_i$ to be the set of all the $(d-k_i)$-wise intersections. It follows that $\abs{\cJ}=\binom{m}{d}$ and $\abs{\cF_i}=\binom{m}{d-k_i}$. For convenience, we set 
    \[C=2^{-\sum_{i=1}^r\bar{w_i}H(e_i)}=\prod_{i=1}^{r}\prod_{e\in c^{-1}(i)}\left(\frac{w(e)}{\bar{w}_i}\right)^{w(e)}.\]
    For each $p^*\in\cJ$, the number of flats in $\cF_i$ that contains $p$ is exactly $\binom{d}{d-k_i}$. 
    By uniform bound, we have
    \[\sum_{i=1}^r\bar{w}_i\left(H(F_{p^*,e_i})-H(e_i)\right)\leq \log_2\left(C\prod_{i=1}^r\binom{d}{d-k_i}^{\bar{w}_i}\right).\]
    Since it is possible to take the random tuple flats $(F_{p^*,e_i})_{i\in [r]}$ such that the marginal distribution of $F_{p^*,e_i}$ is uniform on all the flats in $\cF_i$ that contains $p^*$, the above inequality can be achieved and hence
    \[\eta(p^*)=C\prod_{i=1}^r\binom{d}{d-k_i}^{\bar{w}_i}.\]
    Therefore, the left hand side of the inequality is
    \begin{align*}
        C\binom{m}{d}\prod_{i=1}^r\binom{d}{d-k_i}^{\bar{w}_i}=(1-o(1))C\frac{m^d}{d!}\prod_{i=1}^r\left(\frac{d!}{k_i!(d-k_i)!}\right)^{\bar{w_i}},
    \end{align*}
    and the right hand side is, using that $C_{\cH,w} = C\cdot d!^{\abs{w}-1}\prod_{i=1}^{r}k_i!^{-\bar{w}_i},$
    \begin{align*}
        C_{\cH,w}\prod_{i=1}^r\binom{m}{d-k_i}^{\bar{w}_i}=(1-o(1))Cm^{\sum_{i=1}^r \bar{w}_i(d-k_i)}d!^{\abs{w}-1}\prod_{i=1}^r \left(\frac{1}{k_i!}\right)^{\bar{w}_i}\prod_{i=1}^r\left(\frac{1}{(d-k_i)!}\right)^{\bar{w}_i}.
    \end{align*}
    Note that 
    \[\sum_{i=1}^r \bar{w}_i(d-k_i)=\sum_{e\in E}\sum_{j\in e}w(e)=\sum_{j\in[d]}\sum_{e\ni j}w(e)=d.\]
    Therefore, the two sides only differ by lower order terms as $m$ goes to infinity, and hence $C_{\cH,w}$ is the best possible constant.
\end{proof}

\section{Proof of the key inequality}\label{section:poly}
In this section, we give a proof of \cref{thm:key-ineq}.
To prove the theorem, we first construct related objects that allow us to describe our choice of $b_{p,F}$'s, which is done in \cref{subsec:prior-op}.
Then we define the $b_{p,F}$'s and show that it satisfies the second inequality in the statement of \cref{thm:key-ineq}, which is done in \cref{subsec:joints-of-varieties}, where we also recall the properties we need from Tidor, Yu and Zhao's proof \cite{TYZ22}.
Finally, we show that the $b_{p,F}$'s satisfy the first inequality as well by an application of the Loomis--Whitney inequality and the polynomial method, which is done in \cref{subsec:LW-step}.

Although the properties we need from \cite{TYZ22} are technically from Section 5, the main ideas relevant to this paper are mostly contained in \cite[Section 3]{TYZ22}.
We refer the readers to that section for a more complete discussion of the tools we will use.

\subsection{Derivative operators and linear functionals}\label{subsec:prior-op}
For any injective affine transformation $A:\FF^k\to \FF^d$  and any $\vec{\gamma}\in \ZZ_{\geq 0}^k$, we define the corresponding Hasse derivative $\Hasse^{\vec{\gamma}}_{A}$ by
\[\Hasse^{\vec{\gamma}}_{A}(g) \eqdef \left(\Hasse^{\vec{\gamma}} (g\circ A)\right) \circ A^{-1}.\]
See \cite[Section 4.3]{TYZ22} for a brief introduction of Hasse derivatives in the same context.

We now recall and modify some definitions from \cite{TYZ22}.
Fix a $k$-dimensional flat $F$ throughout this subsection.
As in \cite{TYZ22}, we first fix an arbitrary order of the joints, which we will refer to as the \emph{preassigned order} throughout the paper.
A \emph{handicap} $\vec{\alpha}=(\alpha_p)_{p\in \cJ}$ is an assignment of integers to the joints.
The associated \emph{priority order} given a handicap is a total order on $\cJ\times \ZZ_{\geq 0}$ such that $(p,r)\prec (p',r')$ if and only if either  $r-\alpha_p<r'-\alpha_{p'}$ or $r-\alpha_p=r'-\alpha_{p'}$ and $p$ comes before $p'$ in the preassigned order.

Denote by $(R_F)_{\leq n}$ the space of polynomials on the flat $F$ of degree at most $n$.
For any flat $F$ that contains a point $p$, pick an arbitrary bijective affine transformation $A:\FF^k\to F$ such that $A(0)=p$. 
We set $\totalB^r_{p,F}(n)$ to be the space spanned by the linear functionals from $(R_F)_{\leq n}$ of the form
\[g\mapsto \Hasse^{\vec{\gamma}}_{A} g(p)\]
for some $\vec{\gamma}\in \ZZ_{\geq 0}^{k}$ such that $\abs{\vec{\gamma}}=r$.
It is not hard to see that this matches the definition made in \cite{TYZ22} as $(g\mapsto \Hasse^{\vec{\gamma}}_{A}g)_{\abs{\vec{\gamma}}=r}$ is a basis of the space of derivative operators of order $r$ on $F$.
As a consequence, this definition is independent of $A$.
We prefer this definition as at a later stage of the proof, we would need to specify explicitly a coordinate system using the affine transformation $A$.

Recall that in \cite{TYZ22}, the set $\prioB^r_{p,F}(\vec{\alpha},n)$ is chosen so that for any flat $F$, the disjoint union $\bigcup_{(p',r')\preceq (p,r)}\prioB^{r'}_{p',F}(\vec{\alpha},n)$ is a basis of the vector space $\sum_{(p',r')\preceq (p,r)}\totalB^{r'}_{p',F}(n)$.
Also set $\prioB_{p,F}(\vec{\alpha},n)$ to be the union $\bigcup_{r\in\ZZ_{\geq 0}}\prioB^r_{p,F}(\vec{\alpha},n)$.
As observed there, the size of $\prioB^r_{p,F}(\vec{\alpha},n)$ does not depend on the choice of the basis.
We denote by $B^r_{p,F}(\vec{\alpha},n)$ the size of $\prioB^r_{p,F}(\vec{\alpha},n)$ for our purpose, and set $B_{p,F}(\vec{\alpha},n)$ to be the size of $\prioB_{p,F}(\vec{\alpha},n)$.
In our setting, we will use $B_{p,F}(\vec{\alpha},n)$ instead of $\prioB_{p,F}(\vec{\alpha},n)$ as our argument chooses a basis in a much later part, whereas in \cite{TYZ22}, they make the choice at the very beginning.

The connection of the definitions made here and \cref{thm:key-ineq} is that, roughly speaking, we will choose $\vec{\alpha}$ depending on $n$ so that $b_{p,F} = \lim_{n\to\infty}B_{p,F}(\vec{\alpha},n)/n^{\dim F}$ exists for any $p\in \cJ$ and $F\in \cF_1\cup\cdots\cup \cF_r$ with $p\in F$.
The goal now is to choose $\vec{\alpha}$ appropriately so that the two properties listed in \cref{thm:key-ineq} are met.
The second will follow immediately once we recall the properties from \cite{TYZ22}, whereas the first will be harder to establish.

\subsection{Properties of $B_{p,F}(\vec{\alpha},n)$}\label{subsec:joints-of-varieties}
As $B_{p,F}(\vec{\alpha},n)$ is the size of $\prioB_{p,F}(\vec{\alpha},n)$ for any choice of the basis, the properties of $\abs{\prioB_{p,F}(\vec{\alpha},n)}$ that were proved in \cite{TYZ22} translate to properties of $B_{p,F}(\vec{\alpha},n)$.
We record here several properties that we need. 
For the following properties, let $\cJ$ be the set of joints.

\begin{lemma}[{\cite[Equation (5.1)]{TYZ22}}]\label{lem:sum-of-condition}
    For each flat $F$, any handicap $\vec{\alpha}\in \ZZ_{\geq 0}^{\cJ}$ and any positive integer $n$, we have
    \[\sum_{p\in \cJ\cap F}B_{p,F}(\vec{\alpha},n) =\binom{n+\dim F}{\dim F}.\]
\end{lemma}

\begin{lemma}[{Bounded domain, \cite[Lemma 5.4]{TYZ22}}] \label{lem:bdd-domain}
For each $n\in\NN$ and each flat $F$, there is some $C_{\cJ,F}(n)$ so that if $\vec\alpha \in \ZZ^\cJ$ satisfies $ \alpha_p < \max_{q \in \cJ\cap F} \alpha_q - C_{\cJ,F}(n)$,
then $B_{p,F}(\vec\alpha, n)=0$. 
\end{lemma}

\begin{lemma}[{Monotonicity, \cite[Lemma 5.5]{TYZ22}}] \label{lem:mono}
Let $n$ be a positive integer and $\vec\alpha^{(1)},\vec\alpha^{(2)} \in \ZZ^\cJ$ be two handicaps. 
Suppose $F$ is a flat and $p \in \cJ\cap F$ satisfies 
$\alpha^{(1)}_p-\alpha^{(1)}_{p'}\le \alpha^{(2)}_p-\alpha^{(2)}_{p'}$ for all $p' \in \cJ\cap F$.
Then $B_{p,F}(\vec\alpha^{(1)},n) \leq B_{p,F}(\vec\alpha^{(2)},n)$.
\end{lemma}

\begin{lemma}[{Lipschitz continuity, \cite[Lemma 5.6]{TYZ22}}]\label{lem:lip}
Let $F$ be a flat and $p \in \cJ\cap F$. Suppose we have two handicaps $\vec\alpha^{(1)},\vec\alpha^{(2)} \in \ZZ^\cJ$.
Then
\begin{multline*}
\abs{
B_{p, F}(\vec\alpha^{(1)},n)
- B_{p, F}(\vec\alpha^{(2)},n)
}
\\
\leq
\paren{\binom{n}{\dim F -1} + O_d(n^{\dim F - 2})}\sum_{p'\in \cJ\cap F }\abs{(\alpha^{(1)}_{p'}-\alpha^{(1)}_{p})-(\alpha^{(2)}_{p'}-\alpha^{(2)}_{p})}.	
\end{multline*}
\end{lemma}

These will be the only properties we need to know about the quantities $B_{p,F}(\vec{\alpha},n)$.
As mentioned earlier, we can now easily show that any choice for $b_{p,F} = \lim_{n\to\infty}B_{p,F}(\vec{\alpha},n)/n^{\dim F}$ satisfies the second property in \cref{thm:key-ineq} if the limits exist.
In addition to the second property, we also prove a property that resembles the first property by making an appropriate choice for the $\vec{\alpha}$'s.
For technical reasons, we will instead take the limit only for a subsequence.
As in previous works \cite{YZ23, TYZ22}, we need to temporarily put an extra technical condition.

\begin{definition}
    Let $(\cJ,\cF_1,\ldots, \cF_r)$ be an $\cH$-joints configuration.
    If $F\in \cF_1\sqcup\cdots\sqcup \cF_r$ appears as an entry in $\cT_p$ for some $p\in\cJ$, then we say that \emph{$p$ uses $F$} or \emph{$F$ is used by $p$}.
    The $\cH$-joints configuration is said to be \emph{connected} if for any bipartition $P\sqcup Q = \cJ$ with both $P,Q$ non-empty, there exist $p\in P$, $q\in Q$ and $F\in \cF_1\sqcup\cdots\sqcup \cF_r$ such that $F$ is both used by $p$ and $q$.
\end{definition}

\begin{lemma}\label{lem:equation}
Suppose that $(\cJ,\cF_1,\ldots, \cF_r)$ is a connected $\cH$-joints configuration, and $\sigma:E\to \RR_{> 0}$, $W:J\to \RR_{> 0}$ are any functions.
Then we may choose a sequence of positive integers $(n_s)_{s\in \NN}$ tending to infinity and a sequence of handicaps $(\vec{\alpha}^{(s)})_{s\in \NN}$ so that for each $p\in \cJ$ and $F\in\cF_1\cup\cdots\cup\cF_r$ with $p\in F$, the limit
\[b_{p,F}=\lim_{s\to\infty}\frac{B_{p,F}(\vec{\alpha}^{(s)},n_s)}{n_s^{\dim F}}\]
exists, and the following holds.
\begin{enumerate}
    \item There exists $\lambda\geq 0$ such that for any $p\in \cJ$,
    \[\min_{(F_{p,e})_{e\in E}\in \cT_p}\prod_{e\in E}b_{p,F_{p,e}}^{\sigma(e)}= \lambda\cdot W(p).\]
    \item 
$\sum_{p\in F\cap \cJ}b_{p,F}\leq \frac{1}{(\dim F)!}$
for any $F\in \cF_1\cup\cdots\cup\cF_r$.
\end{enumerate}
\end{lemma}
\begin{proof}
    A large part of this proof is analogous to \cite[Lemma 3.11, Lemma 5.10]{TYZ22}.
    We first choose a handicap $\vec{\alpha}$ for any $n\in \NN$.
    In this proof, by $o(1)$ we mean $o_{\cJ, \cF_1,\ldots, \cF_r, \cH, \sigma, W; n\to\infty}(1)$.
    In the first part of the proof, we will also suppress various dependencies on $n$.
    For any $p\in \cJ$ and $\vec{\alpha}\in \ZZ_{\geq 0}^{\cJ}$, set
    \[W'_p(\vec{\alpha},n) = \frac{1}{W(p)}\min_{(F_{p,e})_{e\in E}\in \cT_p}\prod_{e\in E}\left(\frac{B_{p,F_{p,e}}(\vec{\alpha},n)}{n^{\dim F_{p,e}}}\right)^{\sigma(e)}\]
    and
    \[S_p(\vec{\alpha},n) = \sum_{p\in F}B_{p,F_{p,e}}(\vec{\alpha},n).\]
    When we surpress the dependecies on $n$, we write $W'_p(\vec{\alpha})$ and $S_p(\vec{\alpha})$ instead.
    Observe that $W_p'(\vec{\alpha})$ remains unchanged if we shift each $\alpha_p$ by the same constant, as the priority order induced by $\vec{\alpha}_p$ does not change in this case.
    We also notice that there is a finite number $C$ so that if $\alpha_p<\alpha_{p'}-C$ for some joints $p,p'$ using the same flat $F$, then $B_{p,F}(\vec{\alpha},n)=0$ by \cref{lem:bdd-domain}.
    Therefore $W'_p(\vec{\alpha})=0$ in this case.
    In addition, note that $S_p(\vec{\alpha})$ is a non-negative integer bounded above by $\abs{\cF_1\cup\cdots\cup\cF_r}n^d.$
    What these show is that there are only finitely many possibilities for $((W'_p(\vec{\alpha}),S_p(\vec{\alpha})))_{p\in \cJ}$ for any given $n$.
    We can thus choose $\vec{\alpha}$ so that after the pairs $((W'_p(\vec{\alpha}),S_p(\vec{\alpha})))_{p\in \cJ}$ are sorted in decreasing lexicographical order, the whole sequence of pairs is the largest in lexicographical order.

    Let $J=\abs{\cJ}$, and suppose that 
    \[(W'_{p_1}(\vec{\alpha}),S_{p_1}(\vec{\alpha}))\geq_{\textup{lex}}(W'_{p_2}(\vec{\alpha}),S_{p_2}(\vec{\alpha}))\geq_{\textup{lex}}\cdots\geq_{\textup{lex}}(W'_{p_J}(\vec{\alpha}),S_{p_J}(\vec{\alpha})).\]
    We will now show that $W'_{p_i}(\vec{\alpha})-W'_{p_{i+1}}(\vec{\alpha})\leq \delta_n$ for all $i\in[J-1]$, where $(\delta_n)_{n\in\NN}$ is an appropriately chosen sequence of positive reals that tends to zero sufficiently slowly when $n$ goes to infinity.
    Suppose for the sake of contradiction that this is not the case.
    Then we may take $t$ to be the smallest positive integer such that $W'_{p_t}(\vec{\alpha})-W'_{p_{t+1}}(\vec{\alpha})>\delta_n$.
    Let $\vec{v} = \vec{e}_{p_1}+\cdots +\vec{e}_{p_t}$ and $\vec{\alpha}' = \vec{\alpha}-\vec{v}$.
    Then by the monotonicity of $B_{p,F}(\vec{\alpha},n)$ (\cref{lem:mono}), we have that $B_{p_i,F}(\vec{\alpha}',n)\leq B_{p_i,F}(\vec{\alpha},n)$ for any $i\leq t$, and $B_{p_j,F}(\vec{\alpha},n)\leq B_{p_j,F}(\vec{\alpha}',n)$ for any $j>t$.
    As a consequence, we have $W'_{p_i}(\vec{\alpha}')\leq W'_{p_i}(\vec{\alpha})$ and $S_{p_i}(\vec{\alpha}')\leq S_{p_i}(\vec{\alpha})$ for any $i\leq t$, and $W'_{p_j}(\vec{\alpha}')\geq W'_{p_j}(\vec{\alpha})$ and $S_{p_j}(\vec{\alpha}')\geq S_{p_j}(\vec{\alpha})$ for any $j>t$.

    Now we bound $\abs{W'_p(\vec{\alpha}')-W'_p(\vec{\alpha})}$ for each $p\in \cJ$.
    To do so, note that by \cref{lem:lip}, we have $\abs{B_{p,F}(\vec{\alpha}')-B_{p,F}(\vec{\alpha})}<(\abs{\cJ}+o(1))n^{\dim F-1} = o(n^{\dim F})$ for each $p\in F$.
    Therefore we have $\abs{W'_p(\vec{\alpha}')-W'_p(\vec{\alpha})}=o(1)$ as well.
    Note that in particular, as long as $\delta_n$ tends to zero slowly enough, we see that for any $i\leq t<j$ we still have $W'_{p_i}(\vec{\alpha}')> W'_{p_j}(\vec{\alpha})$.
    Therefore the first $t$ joints remain the first $t$ even if we sort the pairs $((W'_p(\vec{\alpha}'),S_p(\vec{\alpha}')))_{p\in \cJ}$ decreasingly.
    By the monotonicity we showed above, the choice of $\vec{\alpha}$ forces $(W'_{p_i}(\vec{\alpha}'),S_{p_i}(\vec{\alpha}'))=(W'_{p_i}(\vec{\alpha}),S_{p_i}(\vec{\alpha}))$ for all $i\leq t$.
    In particular, we see that $B_{p_i,F}(\vec{\alpha}',n)=B_{p_i,F}(\vec{\alpha},n)$ for every $i\leq t$ and every $F$ containing $p_i$.

    Now suppose that $(W'_{p_j}(\vec{\alpha}'),S_{p_j}(\vec{\alpha}'))\neq(W'_{p_j}(\vec{\alpha}),S_{p_j}(\vec{\alpha}))$ for some $j>t$, then we have either $W'_{p_j}(\vec{\alpha}')>W'_{p_j}(\vec{\alpha})$ or $S_{p_j}(\vec{\alpha}')>S_{p_j}(\vec{\alpha})$.
    In either case, we see that there exists some $F\ni p_j$ such that $B_{p_j,F}(\vec{\alpha}',n)>B_{p_j,F}(\vec{\alpha},n)$.
    By \cref{lem:sum-of-condition}, this means that there exists some $p_i\in F$ such that $B_{p_i,F}(\vec{\alpha}',n)<B_{p_i,F}(\vec{\alpha},n)$, and by the monotonicity established above, we see that $i\leq t$.
    This is a contradiction.
    As a consequence, we see that $(W'_{p_j}(\vec{\alpha}'),S_{p_j}(\vec{\alpha}'))=(W'_{p_j}(\vec{\alpha}),S_{p_j}(\vec{\alpha}))$ for all $j>t$ as well.

    It is now easy to show, by induction, that $(W'_{p}(\vec{\alpha}-m\vec{v}),S_{p}(\vec{\alpha}-m\vec{v}))=(W'_{p}(\vec{\alpha}),S_{p}(\vec{\alpha}))$ for any $p\in \cJ$ and $m\in\NN$.
    However, since the configuration is connected, there exists $i\leq t$ and $j>t$ such that $p_i$ and $p_j$ both use some flat $F$.
    By \cref{lem:bdd-domain}, if $m$ is sufficiently large, then $B_{p_i,F}(\vec{\alpha}-m\vec{v})=0$, showing that $W'_{p_i}(\vec{\alpha})=W'_{p_i}(\vec{\alpha}-m\vec{v})=0$.
    This is a contradiction, and so we indeed have $W'_{p_i}(\vec{\alpha})-W'_{p_{i+1}}(\vec{\alpha})\leq \delta_n$ for every $i\in [J-1]$.
    As a consequence, $\abs{W'_{p}(\vec{\alpha})-W'_{p'}(\vec{\alpha})}\leq J\delta_n = o(1)$ for any two joints $p,p'\in\cJ$.

    Now note that for any $p\in F$ where $p\in\cJ$ and $F\in \cF_1\cup\cdots\cup \cF_r$, we have $B_{p,F}(\vec{\alpha},n)/n^{\dim F}$ is contained in a bounded interval.
    Therefore, we may choose $n_1,n_2,\ldots\in\NN$ tending to infinity so that if $\vec{\alpha}^{(s)}$ is the handicap we chose above corresponding to $n_s$, then
    \[b_{p,F}=\lim_{s\to\infty}\frac{B_{p,F}(\vec{\alpha}^{(s)},n_s)}{n_s^{\dim F}}\]
    exists for any $p\in F$ with $p\in\cJ$ and $F\in\cF_1\cup\cdots\cup\cF_r$.

    Finally, let us show that the $b_{p,F}$'s satisfy the two properties.
    First, for any $p,p'\in\cJ$, we have
    \begin{align*}
        &\abs{\frac{1}{W(p)}\min_{(F_{p,e})_{e\in E}\in \cT_p}\prod_{e\in E}b_{p,F_{p,e}}^{\sigma(e)}-\frac{1}{W(p')}\min_{(F_{p',e})_{e\in E}\in \cT_{p'}}\prod_{e\in E}b_{p',F_{p',e}}^{\sigma(e)}}\\
        =&\lim_{s\to\infty}\abs{W'_p(\vec{\alpha}^{(s)},n_s)-W'_{p'}(\vec{\alpha}^{(s)},n_s)}=0,
    \end{align*}
    showing that 
    \[\lambda\eqdef\frac{1}{W(p)}\min_{(F_{p,e})_{e\in E}\in \cT_p}\prod_{e\in E}b_{p,F_{p,e}}^{\sigma(e)}\geq 0\]
    does not depend on the choice $p\in\cJ$.
    Therefore the first property holds.
    Moreover, by \cref{lem:sum-of-condition}, we have
    \[\sum_{p\in F\cap \cJ}b_{p,F} = \lim_{s\to\infty}\frac{\sum_{p\in F\cap \cJ}B_{p,F}(\vec{\alpha}^{(s)},n_s)}{n_s^{\dim F}} = \lim_{s\to\infty}\frac{\binom{n_s}{\dim F}}{n_s^{\dim F}} = \frac{1}{(\dim F)!}.\]
    Therefore the second property holds as well.    
\end{proof}

We make a brief note that to prove \cref{thm:key-ineq} for connected $\cH$-joints configurations with $w(e)>0$ for all $e\in E$, it now suffices to plug in $\sigma(e) = w(e)/(\abs{w}-1)$ and prove that $\lambda\geq 1$ in \cref{lem:equation}.
We will prove so by showing that 
\[\sum_{p\in\cJ}\min_{(F_{p,e})_{e\in E}\in \cT_p}\prod_{e\in E}b_{p,F_{p,e}}^{\frac{w(e)}{\abs{w}-1}}\geq \frac{1}{d!},\]
and this is where we will use the polynomial method.
In light of this, for each $p\in\cJ$, we will now fix the minimizer $(F_{p,e})_{e\in E}\in \cT_p$, and also fix a choice $A_p$ of an affine transformation that witnesses the joint $p$ corresponding to the minimizer $(F_{p,e})_{e\in E}\in \cT_p$.

\subsection{Lower bound via the polynomial method}\label{subsec:LW-step}
In this subsection, we will conclude the proof of \cref{thm:key-ineq} by showing the desired lower bound using the polynomial method.
Before we begin, we need several notations listed below.

For any set $S$ and any edge $e\in E$, let $\iota^{(e)}:S^{d-\abs{e}}\to S^d$ be the map sending $(s_1,\ldots, s_{d-\abs{e}})$ to $(s_1',\ldots, s_d')$ where $s_i' = 0$ if $i\in e$, and $s_i'=s_{\textup{rank}(i)}$ if $i\not\in e$ where $\textup{rank}(i)$ is the rank of $i$ in the set $[d]\backslash e$.
Here, the rank goes from $1$ to $d-\abs{e}$ starting from the smallest element in $[d]\backslash e$.
We also define $\pi^{(e)}:S^{d}\to S^{d-\abs{e}}$ to be the projection so that $\pi^{(e)}\circ \iota^{(e)}$ is the identity map.
We do not include $S$ in our notation as it will always be clear from context.
Note also that when $S=\FF$, the maps $\iota^{(e)}$ and $\pi^{(e)}$ are linear transformations.

Recall that we have just chosen the tuple $(F_{p,e})_{e\in E}\in \cT_p$ and the affine transformation $A_p$ for a given joint $p$ in the previous subsection.
We are going to choose the basis $\prioB_{p,F_{p,e}}(\vec{\alpha},n)$ according to this affine transformation.
This is different from \cite{TYZ22}, where the choice was made arbitrarily.
Note that for any edge $e\in E$, the affine transformation $A_{p,e}\eqdef A_p\circ\iota^{(e)}$ is a bijective affine transformation $\FF^k\to F_{p,e}$ where $k=\dim F_{p,e} = d-\abs{e}$.
Now we choose $\cG^r_{p,e}(\vec{\alpha},n)\subseteq\ZZ_{\geq 0}^k$ as follows: we list the $\vec{\gamma}\in \ZZ^k_{\geq 0}$ with $\abs{\vec{\gamma}}=r$ in decreasing lexicographical order and record a vector space $W$ that is initially $\sum_{(p',r')\prec (p,r)}\totalB^r_{p',F}(n)$.
Now we go through the list term by term, and we check if the linear functional $\Lambda^{\vec{\gamma}}:g\mapsto \Hasse^{\vec{\gamma}}_{A_{p,e}}g(p)$ is in the vector space $W$.
If not, we put $\vec{\gamma}$ in $\cG^r_{p,e}(\vec{\alpha},n)$ and replace $W$ with the span of $W$ and $\Lambda^{\vec{\gamma}}$.
Note that by the construction, it is clear that when we go over the entire list, we have that 
\[\left\{\Lambda^{\vec{\gamma}}\mid \vec{\gamma}\in \cG^r_{p,e}(\vec{\alpha},n)\right\}\]
is a valid choice of $\prioB^r_{p,F}(\vec{\alpha},n)$ (using our definition of $\prioB^r_{p,F}(\vec{\alpha},n)$ in \cref{subsec:prior-op} with $A=A_{p,e}$).
We set $\cG_{p,e}(\vec{\alpha},n)$ to be the union $\bigcup_{r\in\ZZ_{\geq 0}}\cG^r_{p,e}(\vec{\alpha},n)$.
Observe that $\abs{\cG^r_{p,e}(\vec{\alpha},n)}=B^r_{p,F_{p,e}}(\vec{\alpha},n)$ and $\abs{\cG_{p,e}(\vec{\alpha},n)}=B_{p,F_{p,e}}(\vec{\alpha},n)$ by definition.

We now define $\cG_{p}(\vec{\alpha},n)\subseteq \ZZ^d_{\geq 0}$ from the sets $\left(\cG_{p, e}(\vec{\alpha},n)\right)_{e\in E}.$
Define it to be the subset of $\{\vec{\gamma}\in \ZZ_{\geq 0}^d\mid \abs{\vec{\gamma}}\leq n\}$ containing all $\vec{\gamma}$'s with the property
\[\pi^{(e)}(\vec{\gamma})\in \cG_{p, e}(\vec{\alpha},n)\quad\forall e\in E.\]
The set $\cG_p(\vec{\alpha},n)$ is going to be an important intermediate term in the proof of the lower bound.
We first relate the size of $\cG_p(\vec{\alpha},n)$ to the sizes of $\abs{\cG_{p,e}(\vec{\alpha},n)}$ for $e\in E$.
\begin{lemma}\label{lem:LW-step}
    For any $w:E\to \RR_{\geq 0}$ that covers $\cH$,
    we have
    \[\frac{\abs{\cG_{p}(\vec{\alpha},n)}}{(n+1)^d}\leq \prod_{e\in E}\left(\frac{\abs{\cG_{p,e}(\vec{\alpha},n)}}{(n+1)^{d-\abs{e}}}\right)^{\frac{w(e)}{\abs{w}-1}}.\]
\end{lemma}
\begin{proof}
    Since no edges in $E$ are $d$-uniform, we have that $\abs{w}\geq d/(d-1)>1$.
    Notice that $\pi^{(e)}\left(\cG_{p}(\vec{\alpha},n)\right)\subseteq \cG_{p,e}(\vec{\alpha},n)$ for any $e\in E$ and $\pi_j\left(\cG_{p}(\vec{\alpha},n)\right)\subseteq \{0,\ldots, n\}$ for any $j\in [d]$.
    Apply the Loomis--Whitney inequality to the subset-weight pairs
    \[\left([d]\backslash e, \frac{w(e)}{\abs{w}-1}\right)_{e\in E}\cup\left(\{j\}, \frac{-1+\sum_{e\ni j}w(e)}{\abs{w}-1}\right)_{j\in [d]},\]
    and we get that
    \[\abs{\cG_{p}(\vec{\alpha},n)}\leq \prod_{e\in E}\abs{\cG_{p,e}(\vec{\alpha},n)}^{\frac{\Tilde{w}(e)}{\abs{w}-1}}\prod_{j\in [d]}(n+1)^{\frac{-1+\sum_{e\ni j}w(e)}{\abs{w}-1}}.\]
    The second product can be rewritten as
    \[\left((n+1)^{-d}\prod_{e\in E}(n+1)^{w(e)\abs{e}}\right)^{\frac{1}{\abs{w}-1}}=(n+1)^d\left(\prod_{e\in E}(n+1)^{w(e)\left(-d+\abs{e}\right)}\right)^{\frac{1}{\abs{w}-1}},\]
    and so
    \[\abs{\cG_{p}(\vec{\alpha},n)}\leq(n+1)^d\prod_{e\in E}\left(\frac{\abs{\cG_{p,e}(\vec{\alpha},n)}}{(n+1)^{d-\abs{e}}}\right)^{\frac{w(e)}{\abs{w}-1}}.\]
    We are thus done after rearranging.
\end{proof}

We next lower bound the total size of $\cG_p(\vec{\alpha},n)$ when $p$ runs through all joints.
This step is done via the polynomial method.
\begin{lemma}\label{lem:param-counting}
    \[\sum_{p\in\cJ}\abs{\cG_{p}(\vec{\alpha},n)}\geq \binom{n+d}{d}.\]
\end{lemma}
\begin{proof}
    Let $g$ be any polynomial of degree at most $n$ such that 
    \[\Hasse^{\vec{\gamma}}_{A_p}g(p)=0\quad \forall \vec{\gamma}\in \cG_{p}(\vec{\alpha},n), \, p\in \cJ.\]
    We will show that $g$ has to be the zero polynomial.
    Note that if this is the case, then the inequality follows immediately as the left hand side is the number of linear conditions we put on $g$, and the right hand side is the dimension of the vector space $\FF[x_1,\ldots, x_d]_{\leq n}$.

    Suppose for the contradiction that $g\neq 0$.
    Then there exists some $p\in \cJ$ and $\vec{\gamma}\in \ZZ_{\geq 0}^d$ such that $\Hasse^{\vec{\gamma}}_{A_p}g(p)\neq 0$.
    Choose $p$ and $\vec{\gamma}$ such that $(p,\abs{\vec{\gamma}})$ is the minimum with respect to $\prec$, and we tie break by choosing the $\vec{\gamma}$ with the largest lexicographical order.
    It is clear that $\abs{\vec{\gamma}}\leq n$ but $\vec{\gamma}\not\in \cG_{p}(\vec{\alpha},n)$.
    By the definition of $\cG_{p}(\vec{\alpha},n)$, there exists an edge $e\in E$ such that
    \[\pi^{(e)}(\vec{\gamma})\not\in \cG_{p,e}(\vec{\alpha},n).\]
    For simplicity, set $\vec{\gamma_e} = \pi^{(e)}(\vec{\gamma})$ and let $\Tilde{g}$ be the restriction of $\Hasse^{\vec{\gamma}-\iota^{(e)}(\vec{\gamma_e})}_{A_p}g$ on the flat $F\eqdef \textup{im } A_{p,e}$.
    Then it is clear that $\Hasse^{\vec{\gamma_e}}_{A_{p,e}}\Tilde{g}(p)=\Hasse^{\vec{\gamma}}_{A_p}g(p)\neq 0.$
    Moreover, by the choice of $(p,\vec{\gamma})$, we know that for any $p'\in F$ and any $r\in \ZZ_{\geq 0}$, if $(p',r)\prec (p,\abs{\vec{\gamma}})$, then $\mult(g,p')>r$.
    As we have the lower bound $\mult(\Tilde{g},p')\geq \mult(g,p')-\abs{\vec{\gamma}}+\abs{\vec{\gamma_e}}$, this shows that if $r'\in\ZZ_{\geq 0}$ is such that $(p',r')\prec (p,\abs{\vec{\gamma_e}})$, then $\mult(\Tilde{g},p')>r'$.
    A consequence of this is that for any linear functional $\Lambda\in \sum_{(p',r')\prec (p,\abs{\vec{\gamma_e}})}\totalB_{p',F}(n)$, we have $\Lambda(\Tilde{g})=0$.
    Since $\vec{\gamma_e}\not\in \cG_{p,e}(\vec{\alpha},n)$ and $\Hasse^{\vec{\gamma_e}}_{ A_{p,e}}\Tilde{g}(p)\neq 0$, by the definition of $\cG_{p,e}(\vec{\alpha},n)$, this means that there exists $\vec{\gamma_e}'\in\ZZ_{\geq 0}^{\dim F}$ with a greater lexicographical order than $\vec{\gamma_e}$ such that $\abs{\vec{\gamma_e}'}=\abs{\vec{\gamma_e}}$ and $\Hasse^{\vec{\gamma_e}'}_{A_{p,e}}\Tilde{g}(p)\neq 0$.
    However, setting $\vec{\gamma}' = \vec{\gamma}-\iota^{(e)}(\vec{\gamma_e})+\iota^{(e)}\left(\vec{\gamma_e}'\right)$ gives that $\abs{\vec{\gamma}'} = \abs{\vec{\gamma}}$, $\vec{\gamma}'>_{\textup{lex}}\vec{\gamma}$ and also $\Hasse^{\vec{\gamma}'}_{A_p}g(p)\neq 0$.
    This is a contradiction with the choice of $\vec{\gamma}$.
    Therefore, $g$ must be the zero polynomial, as desired.
    
\end{proof}

With those estimates, we are now ready to prove \cref{thm:key-ineq} using \cref{lem:equation}.

\begin{proof}[Proof of \cref{thm:key-ineq}]
    We first make several cleaning steps.
    Let $\cH'=([d],E')$ be the sub-hypergraph of $\cH$ containing only of edges with $w(e)>0$, and let $\cJ'$ be the joints $p$ in $\cJ$ with $W(p)>0$.
    It is clear that $(\cJ',\cF_1,\ldots, \cF_r)$ is still an $\cH'$-joints configuration, and $w$ still covers $\cH'$.
    We also have
    \[\prod_{e\in E'}b_{p,F_{p,e}}^{\frac{w(e)}{\abs{w}-1}}=\prod_{e\in E}b_{p,F_{p,e}}^{\frac{w(e)}{\abs{w}-1}}.\]
    Moreover, we may split the $\cH'$-joints configuration into connected components and apply \cref{lem:equation} to each of them with $\sigma(e) = w(e)/(\abs{w}-1)$ for every $e\in E'$.
    For each component, we get data $(\vec{\alpha}^{(s)},n_s)_{s\in \NN}$ alongside $b_{p,F}$ for each $p\in \cJ'$, $F\in \cF_1\sqcup\cdots\sqcup \cF_r$ with $p\in F$ and a parameter $\lambda$ from the lemma.
    By setting $b_{p,F}=0$ for any $p\in \cJ\backslash \cJ'$, it suffices to show that $\lambda\geq 1$ for each connected component, so we will assume without loss of generality that the $\cH'$-joints configuration $(\cJ',\cF_1,\ldots, \cF_r)$ is connected to begin with.
    
    For each $p\in \cJ'$, fix $(F_{p,e})_{e\in E'}\in \cT_p$ so that it minimizes the expression $\prod_{e\in E}b^{\frac{w(e)}{\abs{w}-1}}_{p,F_{p,e}}$, and fix $A_p$ to be an affine transformation witnessing the $\cH'$-joint $p$ corresponding to the flats $(F_{p,e})_{e\in E'}\in \cT_p$.
    Then to show that $\lambda\geq 1$, it suffices to show that
    \[\sum_{p\in \cJ'}\prod_{e\in E'}b^{\frac{w(e)}{\abs{w}-1}}_{p,F_{p,e}}\geq \frac{1}{d!}.\]

    By combining \cref{lem:LW-step}, \cref{lem:param-counting} (applying to the $\cH'$-joints configuration $(\cJ',\cF_1,\ldots,\cF_r)$) and noting that $\abs{\cG_{p,e}(\vec{\alpha}^{(s)},n_s)} = B_{p,F_{p,e}}(\vec{\alpha}^{(s)},n_s)$ for each $s\in\NN$, we see that
    \[\sum_{p\in \cJ'}\prod_{e\in E'}\left(\frac{B_{p,F_{p,e}}(\vec{\alpha}^{(s)},n_s)}{(n_s+1)^{\dim F_{p,e}}}\right)^{\frac{w(e)}{\abs{w}-1}}\geq \frac{\binom{n_s+d}{d}}{(n_s+1)^d}.\]
    By taking the limit when $s$ tends to infinity, we get the desired inequality.
\end{proof}

\section{Application on number of graph homomorphisms}\label{section:graphhom}
In this section, we will show how the Friedgut--Kahn theorem (\cref{thm:FK}) and the partial shadow phenomenon (\cref{thm:partial}) follow from the simple $\cH$-joints theorem (\cref{thm:HSimpleJoints}). We will use hyperplanes in general position in $\RR^d$ and their intersections to encode the information of $G$. This gives a generically induced configuration. If we take a generically induced configurations in $\RR^{d+t}$ and project it back to $\RR^d$, we obtain a projected generically induced configuration (as defined in the introduction). 
By applying \cref{thm:HSimpleJoints} to these two types of configurations, we recover the Friedgut--Kahn theorem and show the partial shadow phenomenon, respectively.

\subsection{Friedgut--Kahn}\label{subsec:FK}
In this subsection, we give a proof to the Friedgut--Kahn theorem, which is essentially different from the original one.
\begin{proof}[Proof of \cref{thm:FK} using \cref{thm:HSimpleJoints}]
    Let $G=([m],E')$ be a hypergraph with $\abs{E'}=n$. We label the edges of $\cH$ by $E=\{e_1,\dots,e_r\}$ and color them with $c(e_i)=i$. Since $\cH$ has no isolated vertices, we may pick a weight $w(e)$ that covers $\cH$ with $w(e_1)+\dots+w(e_r)=\rho^*(\cH)$, and let $C_{\cH,w}$ be the constant stated in \cref{thm:HSimpleJoints}. Let $\cM(\cH,G)$ be the set of $d$-subsets of $[m]$ that contains a copy of $\cH$ as a subgraph in $G$.

    We fix $m$ hyperplanes $H_1,\dots,H_m$ in $\RR^d$ in general position. Consider the $\cH$-joints configuration $(\cJ,\cF_1,\dots,\cF_r)$ given by
    \[\cF_i=\{\cap_{j\in e}H_j\mid e\in E',\abs{e}=\abs{e_i}\},\]
    and
    \[\cJ=\{\cap_{j\in A}H_j\mid A\in \cM(\cH,G)\}.\]
    For each $A\in \cM(\cH,G)$, we label the vertices in $A$ so that $A=\{\sigma(1),\dots,\sigma(d)\}$, and $\sigma(e_i)$ is an edge in $G$ for all $i\in [r]$. Note that $p=\cap_{j\in A}H_j$ is indeed an $\cH$-joint, formed by the flats $(\cap_{j\in \sigma(e_i)}H_j)_{i\in [r]}$. This is because the $\cH$-joint is witnessed by the affine transformation that sends $\vec{e_k}$ to the direction of the line $\cap_{j\in [d]\setminus\{k\}}H_j$ for all $k\in [d]$. Thus, by \cref{thm:HSimpleJoints}, we have
    \[M(\cH,G)=\abs{\cJ}\leq C_{\cH,w}\abs{\cF_1}^{w(e_1)}\dots\abs{\cF_r}^{w(e_r)}\leq C_{\cH,w}n^{\rho^*(\cH)}.\qedhere\]
\end{proof}

\subsection{Partial shadow phenomenon}
In this subsection, we  modify the proof slightly to show the partial shadow phenonmenon for general hypergraphs.
\begin{proof}[Proof of \cref{thm:partial} using \cref{thm:HSimpleJoints}]

    To show the theorem, set $E=\{e_1,\dots,e_r\}, c(e_i)=i$, take a covering weight $w$ with $w(e_1)+\dots+w(e_r)=\rho^*(\cH)$, and let $C_{\cH,w}$ be the constant stated in \cref{thm:HSimpleJoints}. Assume $G=([m],E')$ is a hypergraph with $\abs{E'}=n$. Let $\cM(\cC_t(\cH),G)$ be the set of $(d+t)$-subsets of $[m]$ that contains a copy of $\cC_t(\cH)$ as a subgraph in $G$. Also, we denote $e_i^{(t)}=e_i\cup\{d+1,\dots,d+t\}$.

    We may consider the following projected generically induced $\cH$-joints configuration. Let $H_1,\dots,H_m$ be hyperplanes in $\RR^{d+t}$ and let $\pi:\RR^{d+t}\rightarrow\RR^d$ be a generic projection. We take the $\cH$-joints configuration $(\cJ,\cF_1,\dots,\cF_r)$ given by
    \[\cF_i=\{\pi(\cap_{j\in e}H_j)\mid e\in E',\abs{e}=\abs{e_i}+t\},\]
    and
    \[\cJ=\{\pi(\cap_{j\in A}H_j)\mid A\in \cM(\cC_t(\cH),G)\}.\]
    Note that $\cJ$ is a set of points and $\cF_i$ is a set of $(d-\abs{e_i})$-flats in $\RR^d$.
    For each $A\in \cM(\cC_t(\cH),G)$, we label the vertices in $A$ so that $A=\{\sigma(1),\dots,\sigma(d+t)\}$, and $\sigma(e_i^{(t)})$ is an edge in $G$ for all $i\in [r]$. Note that $p=\pi(\cap_{j\in A}H_j)$ is indeed an $\cH$-joint formed by the flats $(\pi(\cap_{j\in \sigma(e_i^{(t)})}H_j))_{i\in [r]}$. Indeed, the $\cH$-joint is witnessed by the affine transformation that sends $\vec{e_k}$ to the direction of the line $\pi(\cap_{j\in [d+t]\setminus\{k\}}H_j)$ for all $k\in [d]$.
    Thus, by \cref{thm:HSimpleJoints}, we have
    \[M(\cC_t(\cH),G)=\abs{\cJ}\leq C_{\cH,w}\abs{\cF_1}^{w(e_1)}\dots\abs{\cF_r}^{w(e_r)}\leq C_{\cH,w}n^{\rho^*(\cH)}.\qedhere\]
\end{proof}

\section{Generalizing generalized H\"older}\label{section:GenHolder}
In this section, we will focus on the special case where the configuration is \emph{axis-parallel}. Namely, every flat is parallel to $\Span(\vec{e_j}\mid j\notin e)$ for some edge $e\in E$. As mentioned in the introduction, when restricting to this type of configurations, \cref{thm:HMultiplicityJoints} and \cref{thm:entropy-key-ineq} becomes the weaker versions of the generalized H\"older's inequality and Shearer's inequality, respectively.

\subsection{Generalized H\"older's inequality}
In this subsection, we prove \cref{thm:GenHolder} using the $\cH$-joints theorem with multiplicities (\cref{thm:HMultiplicityJoints}) and a tensor power trick.
\begin{proof}[Proof of \cref{thm:GenHolder} using \cref{thm:HMultiplicityJoints}]
    Since both sides of the inequality are non-negative and continuous, and the inequality is homogeneous, we may focus on the case where $f_i$ is an integral-valued function for all $i\in [m]$.
    We first prove a weaker statement when $S\subset \RR$.
    For each $i\in [m]$, we take $\cF_i$ to be the multi-set of flats, where we include $f_i(p_i)$ copies of the flat $\pi_i^{-1}(p_i)\subseteq \RR^d$ for each $p_i\in S^{I_i}$. Note that $\abs{\cF_i}=\sum_{p_i\in S^{I_i}}f_i(p_i)$. Thus, we may apply \cref{thm:HMultiplicityJoints} with $\cH=([d],\{I_i\}_{i=1}^m,c)$, where $c(I_i)=i$, and the covering weight $w_1,\dots,w_m$. Note that $c$ is $(\abs{I_1},\dots,\abs{I_m})$-uniform and each flat in $\cF_i$ has dimension $d-\abs{I_i}$. Thus, \cref{thm:HMultiplicityJoints} shows that 
    \[\sum_{p\in\cJ}\eta(p)\leq C_{\cH,w} \prod_{i=1}^m\left(\sum_{p_i\in S^{I_i}}f_i(p_i)\right)^{w_i},\]
    where $\cJ$ is the set of $\cH$-joints, which is also the support of $\prod_{i=1}^m f_i(\pi_i(p))^{w_i}$. Thus, it remains to show that $\eta(p)\geq \prod_{i=1}^m f_i(\pi_i(p))^{w_i}$. First, note that $H(e_i)=0$ since there is only one edge for each color. Also, if we pick any $F_i\in \cF_i$ with $p\in F_i$ for each $i\in [m]$, $F_1,\dots,F_m$ form an $\cH$-joint at $p$. Thus, for each fixed $p\in\cJ$, we may take $F_i$ to be a random flat from $\cF_i$ which contains $p$, uniformly and independently for each $i$. It follows that \[\log_2\eta(p)\geq \sum_{i=1}^m w_iH(F_i)=\sum_{i=1}^m w_i\log_2 f_i(\pi_i(p)),\]
    and we get that
    \begin{align}\label{eq:weaker-Holder}        
    \sum_{p\in S^d}\prod_{i=1}^{m}f_i(\pi_i(p))^{w_i}\leq C_{\cH,w} \prod_{i=1}^m\left(\sum_{p_i\in S^{I_i}}f_i(p_i)\right)^{w_i}
    \end{align}
    as long as $S\subseteq \RR$.
    The same inequality also holds even if $S$ is not a subset of $\RR$, as we can map $S$ injectively into $\RR$ and apply the inequality.

    To remove the extra factor $C_{\cH,w}$, for each $n\in\NN$ we define $f^{(n)}_i:(S^n)^{I_i}\to\RR_{\geq 0}$, the $n$-th tensor power of $f_i$, as follows.
    For each $q_i\in (S^n)^{I_i}=(S^{I_i})^n$, write $q_i=(q_i^{(1)},\ldots, q_i^{(n)})$ where $q_i^{(1)},\ldots, q_i^{(n)}\in S^{I_i}$.
    Then set $f^{(n)}_i(q_i) = f_i(q_i^{(1)})\cdots f_i(q_i^{(n)})$.
    Applying \cref{eq:weaker-Holder} with $S^n$ in place of $S$ and $f_i^{(n)}$ in place of $f_i$ for each $i\in [m]$, we get
    \begin{align}\label{eq:tensored-Holder}
        \sum_{q\in (S^n)^d}\prod_{i=1}^{m}f_i^{(n)}(\pi_i(q))^{w_i}\leq C_{\cH,w} \prod_{i=1}^m\left(\sum_{q_i\in (S^n)^{I_i}}f_i^{(n)}(q_i)\right)^{w_i}.
    \end{align}
    Note that
    \[\sum_{q_i\in (S^n)^{I_i}} f_i^{(n)}(q_i)= \sum_{q_i^{(1)},\ldots, q_i^{(n)}\in S^{I_i}}f_i(q_i^{(1)})\cdots f_i(q_i^{(n)}) = \left(\sum_{p_i\in S^{I_i}}f_i(p_i)\right)^n,\]
    so the right hand side of \cref{eq:tensored-Holder} is 
    \[C_{\cH,w} \left[\prod_{i=1}^m\left(\sum_{p_i\in S^{I_i}}f_i(q_i)\right)^{w_i}\right]^n.\]
    Similarly, if we write $q\in (S^n)^d=(S^d)^n$ as $(q^{(1)},\ldots, q^{(n)})$ where $q^{(1)},\ldots, q^{(n)}\in S^d$, then
    \[\sum_{q\in (S^n)^d}\prod_{i=1}^{m}f_i^{(n)}(\pi_i(q))^{w_i} = \sum_{q^{(1)},\ldots, q^{(n)}\in S^d}\prod_{i=1}^{m}\prod_{j=1}^{n}f_i(\pi_i(q^{(j)}))^{w_i}=\left(\sum_{p\in S^d}\prod_{i=1}^{m}f_i(\pi_i(p))^{w_i}\right)^n.\]
    Therefore by \cref{eq:tensored-Holder},
    \[\left(\sum_{p\in S^d}\prod_{i=1}^{m}f_i(\pi_i(p))^{w_i}\right)^n\leq C_{\cH,w} \left[\prod_{i=1}^m\left(\sum_{p_i\in S^{I_i}}f_i(q_i)\right)^{w_i}\right]^n,\]
    and by taking the $n$-th root of both sides when $n$ goes to infinity, we get the desired inequality.
\end{proof}

\subsection{Shearer's inequality}\label{subsec:shearer}
In this subsection, we prove Shearer's inequality (\cref{thm:Shearer}) using the geometric Shearer's inequality (\cref{thm:entropy-key-ineq}) and a tensor power trick.
\begin{proof}[Proof of \cref{thm:Shearer} using \cref{thm:entropy-key-ineq}]
We first prove a weaker inequality when $X_1,\dots,X_d$ are $\RR$-valued random variables with finite support. Similar as before, we may consider the hypergraph $\cH=([d],\{I_i\}_{i=1}^m,c)$, where $c(I_i)=i$, and the covering weight $w_1,\dots,w_m$. Let $p=(X_1,\dots,X_d)$ be a random point and let $F_{p,I_i}$ be the $(d-\abs{I_i})$-flat that agrees with $X_{I_i}$ in the coordinates indexed by $I_i$. It follows that $(F_{p,I_i})_{i\in [r]}\in\cT_p$. Note that, since there is only one edge for each color, $F_{p,e_i}=F_{p,I_i}$ and $H(e_i)=0$. Moreover, $H(F_{p,e_i}\mid p)=0$ since $F_{p,e_i}$ is uniquely determined by $p$. Therefore, \cref{thm:entropy-key-ineq} implies
\begin{align}\label{eq:weaker-Shearer}
    H(X_1,\dots,X_d)\leq \sum_{i=1}^r w_iH(X_{I_i})+\log_2 C_{\cH,w}.
\end{align}
    
For random variables with finite support that are not $\RR$-valued, we may compose them with an injection from their supports to $\RR$ and apply the argument above. Note that entropy is preserved under injection. Thus, \cref{eq:weaker-Shearer} holds for any random variables with finite support.

    In order to remove the term $\log_2 C_{\cH,w}$, we apply the following tensor power trick. We resample $(X_1,\dots,X_d)$ independently $n$ times and get $(X^{(j)}_1,\dots,X^{(j)}_d)$ for each $j\in [n]$. By applying \cref{eq:weaker-Shearer} with $(X^{(1)}_i,\dots,X^{(n)}_i)$ in place of $X_i$, we get
    \[nH(X_1,\dots,X_d)\leq n\sum_{i=1}^r w_iH(X_{I_i})+\log_2 C_{\cH,w}.\]
    Dividing by $n$ on both sides and taking $n$ goes to infinity, we can conclude that
    \[H(X_1,\dots,X_d)\leq \sum_{i=1}^r w_iH(X_{I_i}).\qedhere\]
\end{proof}

\subsection{Removing the extra factors without tensor power tricks}
In both the generalized H\"older's inequality and Shearer's inequality, the direct applications of our theorems give extra factors. 
We removed the factors by tensor power tricks in our proofs above. 
Here, we briefly remark that there is another way to prove these two results without using the trick, which is done by modifying the polynomial method part. 
As it requires some effort, we do not go through all the details of this. 
The key idea is to replace the space of polynomials of total degree at most $n$ by the space of polynomials where the degree of each variable $x_i$ is at most $n$ for each $i\in [d]$. 
With this modification, we may obtain a version of \cref{thm:key-ineq} for axis-parallel configurations with $\frac{1}{d!}$ and $\frac{1}{k_i!}$ being replaced by $1$ in the theorem statement once we modify the definitions and lemmas in \cref{section:poly} properly. 
The remaining deduction can be done in the same way.

\section*{Acknowledgement}
We would like to thank Zeev Dvir, Noga Alon and Jonathan Tidor for helpful comments on the earlier drafts of the paper.

\bibliographystyle{amsplain0}
\bibliography{ref_joints}
\end{document}